\documentclass[a4paper,12pt]{article}
\usepackage{dsfont}
\usepackage{amsmath}
\usepackage{enumitem}
\usepackage{amssymb,amsthm}
\usepackage{mathrsfs,amsfonts}
\usepackage{color}
\usepackage{graphicx}
\usepackage{bbm}
\usepackage{ulem}
\usepackage{hyperref}
\usepackage[capitalise]{cleveref}

\allowdisplaybreaks

\allowdisplaybreaks

\newtheorem{con0}{Theorem}[section]
\newtheorem{thm0}{Theorem}[section]
\newtheorem{exa0}{Theorem}[section]

\newtheorem{con1}[con0]{Condition}
\newtheorem{def1}[thm0]{Definition}
\newtheorem{lem1}[thm0]{Lemma}
\newtheorem{thm1}[thm0]{Theorem}
\newtheorem{cor1}[thm0]{Corollary}
\newtheorem{pro1}[thm0]{Proposition}
\newtheorem{rem1}[thm0]{Remark}
\newtheorem{ass1}[thm0]{Assumption}
\newtheorem{exa1}[exa0]{\it{Example}}

\def\bglemma{\begin{lem1}}\def\edlemma{\end{lem1}}

\def\bgproposition{\begin{pro1}}\def\edproposition{\end{pro1}}

\def\benumerate{\begin{enumerate}}\def\eenumerate{\end{enumerate}}
\def\bitemize{\begin{itemize}}\def\eitemize{\end{itemize}}

\def\beqlb{\begin{eqnarray}}\def\eeqlb{\end{eqnarray}}
\def\beqnn{\begin{eqnarray*}}\def\eeqnn{\end{eqnarray*}}

\def\eqref#1{{\rm(\ref{#1})}}

\def\ar{\!\!\!&}

\def\proof{\noindent{\it Proof.~}}
\def\qed{\hfill$\square$\smallskip}

\def\mrm{\mathrm}\def\mbb{\mathbf}\def\mcr{\mathscr}
\def\mbb{\mathbb}\def\mds{\mathds}

\def\d{\mrm{d}}\def\e{\mrm{e}}

\def\I{\mds{1}}

\newcommand{\R}{\mathbb R}

\begin{document}

	\centerline{\Large\bf A localized coupling approach to}
	
	\smallskip
	
	\centerline{\Large\bf interacting continuous-state branching processes}

	\bigskip
	
	\centerline{Shukai Chen$^1$, Pei-Sen Li$^2$ and  Jian Wang$^3$}
	
	\medskip
	
	{\small\it

       \centerline{$^1$School of Mathematics and Statistics, Fujian Normal University, Fuzhou 350007, China,}

       \centerline {\tt skchen@fjnu.edu.cn}

		\smallskip
        \centerline{$^2$School of Mathematics and Statistics, Beijing Institute of Technology, Beijing}
		
		\centerline{100872, China, {\tt peisenli@bit.edu.cn}}
		
		\smallskip
        \centerline{$^3$School of Mathematics and Statistics \& Key Laboratory of Analytical Mathematics}

        \centerline{and Applications (Ministry of Education) \& Fujian Provincial Key Laboratory}

        \centerline{
of Statistics and Artificial Intelligence, Fujian Normal University, Fuzhou 350007, China,}

         \centerline{\tt jianwang@fjnu.edu.cn}
		
	}
	
	\bigskip
	
	{\narrower{\narrower
			
			\noindent{\textit{Abstract:}} We introduce
a class of
continuous-state branching processes with immigration, predation and competition, which can be viewed as a combination of the classical Lotka-Volterra model and continuous-state branching processes with competition that 
were introduced by Berestycki, Fittipaldi, and Fontbona (Probab. Theory Relat. Fields, 2018). This model
can be constructed as 
a unique strong solution to a
class of
two-dimensional stochastic differential equations with jumps. We establish sharp conditions for the uniform ergodicity in
the total variation of this model. Our proof relies on a novel, localized Markovian coupling approach, which
is
of its own interest in the ergodicity theory of Markov processes
with interactions.
			
			\bigskip
			
			
			\smallskip
			
			\noindent{\textit{MSC {\rm(2020)} Subject Classification:}} 60J80, 60J25, 60G51, 60G52
			
			\par}\par}

	\bigskip

\section{Introduction and main results}

 \setcounter{equation}{0}
\subsection{Background and motivation}
In this work, we investigate
the uniform ergodicity of
a class of stochastic models integrating the Lotka-Volterra system with
continuous-state branching processes. We begin with a brief introduction to both components
in these models.

The Lotka-Volterra system of ordinary differential equations describes the biological interactions between a predator species and its prey species. It was independently developed by  Lotka \cite{Lotka20} and Volterra \cite{Volterra26},  and  has become a foundational model in theoretical ecology. The system is defined by the following pair of equations:
\beqlb\label{determin LV}
\left\{\begin{array}{ll}
    \displaystyle \frac{\d X_t}{\d t} = X_t(a-bY_t), \medskip \\
    \displaystyle \frac{\d Y_t}{\d t} = Y_t( -c + dX_t),
\end{array}\right.
\eeqlb
where $(X_t)_{t\ge0}$ and $(Y_t)_{t\ge0}$ represent the prey and predator population densities, respectively, and $a, b, c, d$ are positive constants. Here, $a$ denotes the intrinsic growth rate of the prey population in the absence of predators, while $c$ represents the mortality rate of the predator population in the absence of prey. The terms $-bX_tY_t$ and $dX_tY_t$ model the loss of prey due to predation and the corresponding gain in predators resulting from prey consumption, respectively.

To model the effects of environmental factors, Arnold \cite{Arnold1973} proposed the stochastic Lotka-Volterra system by introducing white noise perturbations into coefficients
$aX_t$ and $-cY_t$ of \eqref{determin LV}. This  system is constructed as the solution to the following stochastic differential equation (SDE)
\begin{equation}\label{slv}
	\begin{cases}
		\d X_t = X_t(a - bY_t)\,\d t + \sigma_1 X_t \,\d W^{(1)}_t, \\
		\d Y_t = Y_t(
-c+dX_t)\,\d t + \sigma_2 Y_t \,\d W^{(2)}_t,
	\end{cases}
\end{equation}
  where $(W^{(1)}_t, W^{(2)}_t )_{t\ge0}$ denotes a two-dimensional Brownian motion with $ \sigma_i > 0,$ $i=1,2,$ representing the noise intensity.  Probabilists have increasingly studied how the random noise alters
  long-term behaviors (including extinction, coexistence, ergodicity) of  stochastic models compared to their deterministic counterparts. For related works, we refer to Bao et al. \cite{BMYY11}, Hening and Nguyen \cite{HN18}, Mao \cite{M11}, Nguyen and Yin \cite{NY17}, Rudnicki \cite{Rud03}, Zhu and Yin \cite{ZY09} and references therein.

Alternatively, in this work, instead of considering \eqref{slv}, we focus on a different type of stochastic perturbation. Specifically, we consider a model that combines the Lotka-Volterra system with
continuous-state branching processes.
We next recall some background and related work on branching processes.
The {\it Galton-Watson branching process}
was proposed in Galton and Watson \cite{GW74} to model  the extinction of family names. Perhaps it is the simplest stochastic model
for a population evolving in time. The basic assumption is that
each individual independently produces offspring according to an identical probability distribution.
More precisely, the model can be described as follows:
\beqnn
Z(n) = \sum_{i=1}^{Z(n-1) }\xi_{n,i}, \quad n\ge 1,
\eeqnn
where $(\xi_{n,i})_{n,i\ge1}$
is a sequence of i.i.d. random variables taking values in $\mbb{N}. $ Similar to the relation between random walk and Brownian motion, it is natural
to consider branching models with continuous-state space $\mbb{R}_+$. Feller \cite{Fel51} first demonstrated that, under a moment condition on the offspring distribution, the scaling limit of  Galton-Watson branching processes converges to a diffusion process solving the following
stochastic  equation:
\begin{equation}\label{Feller}
	X_t = x- b\int_0^tX_s\,\d s+\int_0^t\sqrt{2cX_s} \,\d
	B_s,\end{equation}
which is now referred to as the {\it Feller
branching diffusion}.
Note that the diffusion coefficient in \eqref{Feller} takes the form of a square root, differing from the linear diffusion coefficient in \eqref{slv}.
The fundamental property of the solution
to \eqref{Feller} is that its transition semigroup $(Q_t)_{t\ge0}$ satisfies the  branching property, i.e., for any $t\geq 0$ and $x,y\ge0,$
\beqlb\label{bp}
Q_t(x, \cdot)*Q_t(y, \cdot)=Q_t(x+y, \cdot),
\eeqlb
where ``\,$*$'' denotes the convolution operator.
This property means that different particles in population act independently of each other. In general, a c\'adl\'ag,  Markov process defined on $\mbb{R}_+$ is called a {\it continuous-state branching process} (CB-process for short) if its transition semigroup satisfies \eqref{bp}.
The construction of CB-processes with immigration mechanism via stochastic integral equations was given by Dawson and Li \cite{DaL06}. Their work generalizes Feller's equation \eqref{Feller} through the incorporation of a jump component:
\begin{equation}\label{DaL sde}
	X_t = x- b\int_0^tX_s\,\d s+\int_0^t\sqrt{2cX_s}\,\d
	B_s+\int_0^t\int_0^\infty\int_0^{X_{s-}}\xi\,\tilde{M}(\d s,\d \xi,\d u)+\eta_t,
\end{equation}
where $M(\d s,\d z,\d u)$ is a Poisson random measure on $(0,\infty)^3$ with intensity $\d s n(\d \xi)\d u$, $\tilde{M}(\d s,\d \xi,\d u)=M(\d s,\d z,\d u)-\,\d s\, n(\d \xi)\,\d u$ is the compensated measure, and $(\eta_t)_{t\ge0}$ is a subordinator defined by
$$
\eta_t=\gamma t+\int_0^t\int_0^\infty \xi\,N(\d s,\d \xi)
$$
that constitutes the immigration component. The immigration mechanism is naturally motivated by ergodicity considerations, as it ensures that the stationary distribution is not
degenerate at zero. Here and in the sequel, we make the convention that
$$
\int_y^x=\int_{(y,x]}\quad \text{and}\quad \int_x^\infty=\int_{(x,\infty)},\quad x\ge y\in\mbb{R}.
$$
One could replace the diffusion term in \eqref{DaL sde} by the stochastic integral $\sqrt{2c}\int_0^t\int_0^{X_s}\,W(\d s,\d u)$
 in terms of a white noise $\{W(\d s,\d u)\}$ as in Dawson and Li \cite{DaL12}.  This resulting equation not only describes an equivalent process for any fixed $x \ge 0$, but also defines an equivalent flow.

Later, Berestycki et al.\ \cite{BFF18} incorporated a competition structure into CB-processes. Specifically, they constructed the \textit{continuous-state branching process with competition} (CBC-process for short) as the unique strong solution to
\beqlb\label{CBC}
\begin{split}
Y_t  =& y+ \sqrt{2c}\int_0^t\int_0^{Y_s}\, W(\d s,\d u)  - \int_0^t [bY_s+g(Y_s)]\,\d s\\
&+ \int_0^t\int_0^1 \int_0^{Y_{s-}} \xi\, \tilde{M}(\d s,\d \xi,\d u) + \int_0^t\int_1^\infty \int_0^{Y_{s-}} \xi \,M(\d s,\d \xi,\d u),\end{split}
\eeqlb
which generalized the SDE construction for the CB-process established by Dawson and Li \cite{DaL12}.
Here, $y\in \mbb{R}_+$, $b\in\mathbb{R}$ and $g$ is a continuous and non-decreasing function with $g(0)=0$, describing the intensity of competition between individuals.
  In particular, when $g(x)=ax^2$ for some $a>0$, the solution corresponds to the logistic branching process introduced by Lambert \cite{Lam05}, where the quadratic term quantifies the negative interactions between each pair of individuals within the population. When $g(x)\equiv0$, the equation \eqref{CBC} reduces to the equation for the CB-process established in Dawson and Li \cite{DaL12}, whose transition semigroup satisfies \eqref{bp}.

The system studied in this paper is formally defined as follows. Consider a filtered probability space $(\Omega, \mathcal{F}, (\mathcal{F}_t)_{t\geq 0}, \mbb{P})$ satisfying the usual conditions. Let $\{W_1(\d s,\d u)\}$ and  $\{W_2(\d s,\d u)\}$ be two independent $(\mathcal{F}_t)$-adapted white noises on $(0,\infty)^2$, both with the same intensity measure  $\d s\,\d u$. Let $\{N_1(\d s, \d \xi, \d u)\}$ and $\{N_2(\d s, \d \xi, \d u)\}$ be two independent $(\mathcal{F}_t)$-Poisson random measures on $(0,\infty)^3$ with intensity measures $\d s\,n_1(\d \xi)\,\d u$ and $\d s\,n_2(\d \xi)\,\d u$, respectively. We assume that they are also independent of $\{W_1(\d s,\d u)\}$ and  $\{W_2(\d s,\d u)\}$. Here, $n_i(\d \xi)$ ($i=1,2$) are $\sigma$-finite measures on $(0,\infty)$ satisfying
\beqlb\label{measure moment}
\int_{(0,\infty)} (\xi \wedge \xi^2) \, n_i(\d \xi) < \infty.
\eeqlb
$\tilde{N}_i(\d s, \d \xi, \d u)$ $(i=1,2)$ denote their compensated counterparts; i.e., $\tilde{N}_i(\d s, \d \xi, \d u)=N_i(\d s,\d \xi,\d u)-\d s \,n_i(\d \xi)\,\d u$.   A \textit{continuous-state branching
process with immigration, predation and competition} (CBIPC-process for short)  is defined by the following system of stochastic equations with any ${\mathcal{F}}_0$-measurable random vector $(X_0, Y_0) \in \mathbb{R}_+^2$:
\beqlb\label{main sde}
\left\{\begin{aligned}
X_t &= X_0 + \int_0^t \left(-b_1 X^{\alpha_1}_s + a_1 X_s + \gamma_1\right)\, \d s + \sqrt{2\sigma_1}\int_0^t\int_0  ^{X_s}\,W_1(\d s,\d u) \\
&\quad + \int_0^t \int_0^\infty \int_0^{X_{s-}} \xi  \,\tilde{N}_1(\d s, \d \xi, \d u), \\[2mm]
Y_t &= Y_0 + \int_0^t \left(k X_sY_s - b_2 Y^{\alpha_2}_s + a_2 Y_s + \gamma_2\right)\, \mathrm{d}s \\
&\quad + \sqrt{2\sigma_2}\int_0^t\int_0^{Y_s}\, W_2(\d s,\d u) + \int_0^t \int_0^\infty \int_0^{Y_{s-}} \xi  \,\tilde{N}_2(\d s, \d \xi, \d u),
\end{aligned}\right.
\eeqlb where $b_i>0, a_i\in\mbb{R}, \alpha_i>0, \gamma_i>0, \sigma_i\ge0$ for $i=1,2$, and $k\in\mbb{R}$.

\bgproposition\label{unique strong solution}
There exists a unique nonnegative strong solution $(X_t,Y_t)_{t\ge0}$ to \eqref{main sde}.
\edproposition

The proof of Proposition \ref{unique strong solution} will be presented in the Appendix. For the existence and the uniqueness of the strong solution to the SDE
\eqref{main sde}, to our knowledge, only a few references exist for related results.
Under suitable Lipschitz conditions, the analysis of solutions to such kind
stochastic equations involving more complex mechanisms was demonstrated in Ren et al. \cite{RXYZ22}, where extinction properties are also established. Under symmetry assumptions, solutions to these
stochastic equations with vanishing jump terms and immigration components have been investigated in M\'el\'eard and coauthors \cite{CM10,MV12}. The quasi-stationary distribution of the process has also been studied therein.

\subsection{Main result}
Let $\mbb{R}_+ = [0,\infty)$. Denote by $\mcr{P}(\mbb{R}^2_+)$ the set of Borel probability measures on $\mbb{R}^2_+$, and let  $\|\cdot\|_{\rm Var}$ be the total variation norm on $\mcr{P}(\mbb{R}^2_+)$.  A key observation is that the total variation norm corresponds to the Wasserstein metric induced by the discrete metric $\I_{\{z \neq \tilde{z}\}}$ 
for all $z,\tilde z\in \mbb{R}^2_+$. Specifically, for any $\gamma, \eta \in \mcr{P}(\mathbb{R}^2_+)$, we have the identity:
$$
\|\gamma-\eta\|_{\rm Var}
 =
2\inf_{\pi\in \mcr{C}(\gamma,\eta)}\int_{\mbb{R}_+^2\times\mbb{R}_+^2} \I_{\{z \neq \tilde{z}\}} \,\pi(\d z,\d \tilde{z}),
$$
where $\mcr{C}(\gamma,\eta)$ is the collection of all probability measures on $\mbb{R}^4_+$ with marginals $\gamma$ and $\eta$. We say that a CBIPC-process $(X_t,Y_t)_{t\ge0}$ or its transition semigroup $(P_t)_{t\ge 0}$ is \textit{uniformly ergodic} in the  total variation distance with \textit{rate} $\lambda_*>0$, if there exist a unique stationary probability distribution $\gamma$ and a constant $C>0$ such that
 \beqlb\label{eq1.10}
\|\eta P_t-\gamma\|_{\rm Var}\le C\e^{-\lambda_*t},
 \quad
t\ge 0, \eta\in \mcr{P}(\mbb{R}^2_+).
 \eeqlb
    The uniform ergodicity \eqref{eq1.10} follows by standard arguments (see also the proof of Theorem \ref{main result1} in Section \ref{sect:3} for more details) if there exists a constant $C> 0$ such that
$$
\sup_{x,\tilde{x},y,\tilde{y}}\|\delta_{(x,y)} P_t-\delta_{(\tilde{x},\tilde{y})} P_t\|_{\rm Var}  \le C\e^{-\lambda_*t}, \quad t\ge 0,
$$
where $\delta_{(x,y)}$ and $\delta_{(\tilde{x},\tilde{y})}$ denote the Dirac measures concentrated on $(x,y)$ and $(\tilde{x},\tilde{y})$, respectively.

To state our contribution, we need the following condition:

\begin{con1}\label{Condition}
Assume that the following two conditions hold:

\begin{itemize}

\item [\rm(i)] $\alpha_i>1, i=1,2$.

\item [\rm(ii)] For each $i=1,2$, either $\sigma_i>0$ or there are constants $\theta_i\in(1,2)$ and $C_i>0$ such that
\beqlb\label{jump cond}
\int_0^x \xi^2\,n_i(\d \xi)\ge C_ix^{2-\theta_i},\quad x\in (0,1].
\eeqlb

\end{itemize}
\end{con1}

The main theorem in this paper is as follows.

\begin{thm1}\label{main result}
 Suppose that Condition $\ref{Condition}$ holds. Then the CBIPC-process $(X_t,Y_t)_{t\ge0}$ to \eqref{main sde} is uniformly ergodic in the  total variation distance.
\end{thm1}

\begin{rem1}\rm
We make some comments on Theorem \ref{main result}.
\begin{itemize}
\item[{\rm(i)}] We argue that Condition \ref{Condition}(i) is, in a certain sense, sharp for the uniform ergodicity. To illustrate this, consider a two-dimensional Cox-Ingersoll-Ross model, defined as the pathwise unique nonnegative solution $(Z^1_t, Z^2_t)_{t \ge 0}$ to the following system of stochastic equations with initial value $(z_1, z_2) \in \mathbb{R}_+^2$:
 \beqnn
\left\{\begin{array}{ll}
 \displaystyle Z^1_t=z_1+\int_0^t(-b_1Z^1_s+\gamma_1)\,\d s+ \sqrt{2\sigma_1}\int_0^t\int_0^{Z^1_s}\, W_1(\d s,\d u),\cr
 \displaystyle Z^2_t=z_2+\int_0^t(-b_2Z^2_s+\gamma_2)\,\d s
 +\sqrt{2\sigma_2}\int_0^t\int_0^ {Z^2_s} \,W_2(\d s,\d u),
\end{array}\right.
 \eeqnn
where the parameters and the noises are the same as those in \eqref{main sde}. Define $\zeta_{1}=\inf\{t\ge0: Z^1_t=1\}$. According to Duhale et al. \cite[Corollary 9]{DFM14}, for any $z_1>1$,
$$
\mbb{E}_{z_1}[\zeta_1]=\int_0^\infty\frac{\mrm{e}^{-\xi}-\mrm{e}^{-z_1\xi}}{b_1\xi+\sigma_1\xi^2}
\exp\left(\int_0^\xi\frac{\gamma_1}{b_1+\sigma_1u}\,\d u\right)\,\d\xi.
$$
Letting $z_1\rightarrow\infty$, we conclude that $\sup_{z_1>1}\mbb{E}_{z_1}[\zeta_1]=\infty$.
 This, together with  Mao \cite[Lemma 2.1]{MaY02}  implies
 that $(Z^1_t)_{t\ge0}$ is not uniformly ergodic, and hence neither is $(Z^1_t, Z^2_t)_{t \ge 0}$.

\item[{\rm(ii)}]  As a preliminary effort in the study of the uniform ergodicity for interacting branching systems, the interaction structure in \eqref{main sde} is comparatively simpler than that of \eqref{determin LV}. More complex interaction structures will be addressed separately. Nevertheless, it should be emphasized that unidirectional interactions are also of significant importance in the contexts of mathematical biology and evolutionary dynamics. For example, Bovier and Hartung \cite{BH23} analyzed a two-component system of coupled Fisher-KPP equations, modeling the evolution of two population types  $(N^A_t,N^B_t)_{t\ge0}$. By defining the total population mass $N^T=N^A+N^B$, the process $(N^T_t,N^A_t)_{t\ge0}$ is transformed into a unidirectionally interacting system. Several related studies, such as Chen et al. \cite{CTW17} and Holzer and Scheel \cite{HS13}, also employ reduced-dimensional representations to simplify interactions.
\end{itemize}
\end{rem1}

\subsection{Approach and novelties}
The ergodicity of the Lotka-Volterra system perturbed by linear Brownian noise \eqref{slv} has been extensively studied. The existing literature, including works by Bena\"{\i}m et al. \cite{BHS08} and Hening and coauthors \cite{HN18, HNC21, HNTU25}, relies on the Foster-Lyapunov criteria developed by Meyn and Tweedie \cite{MT93}. A key step in this approach involves verifying the irreducibility of skeleton chains.  However, this verification becomes particularly challenging in models with jumps, which are not covered by the processes considered in the aforementioned works.

An alternative approach that circumvents the need for irreducibility is the probabilistic 
coupling method. Rather than analyzing the skeleton chain, this method constructs two copies of the process and studies their meeting time, thereby providing direct estimates of convergence in the Wasserstein-type metric. For a comprehensive introduction to this method, we refer to Chen \cite{C05}, Lindvall \cite{Lin92} and the references therein. Recently, this method has also been widely employed to study the ergodicity of SDE with jumps; see, e.g., Li et al. \cite{LLWZ25}, Li and Wang \cite{LW20}, Luo and Wang \cite{LW19} and Wang \cite{W11}. For this reason, we aim to adopt a probabilistic coupling argument in our work to investigate the ergodicity of the models under consideration.

However, as will be explained in the following Remark \ref{rem for drift term}, the existing coupling framework is not directly applicable to our setting. This limitation motivates us to develop a novel approach, which we term the {\it localized coupling method}.

\begin{rem1}\label{rem for drift term}\rm
\begin{itemize}
\item[{\rm(a)}] Due to the mutual independence of $N_1$ and $N_2$, establishing the ergodicity for such high-dimensional jump SDEs remains a largely open problem, as it invalidates standard coupling methods. Under the assumption that the marginal distributions of the Poisson random measures are correlated, Wang and coauthors \cite{LMW21,LW19,ScW12}, Majka \cite{Maj17} and F.Y. Wang \cite{W11} derived several results by using the coupling approach. To clarify the distinction between these references and our setting, consider the $d$-dimensional $(d\ge1)$ jump SDE studied in Luo and Wang \cite{LW19}:
\begin{align}\label{OU-proc}
&\d Z_t=b(Z_t)\,\d t+\d L_t,\quad Z_0\in\mbb{R}^d,
\end{align}
where
$$
L_t=\int_0^t\int_{|\xi|\le1}\xi\,\tilde{N}(\d s,\d\xi)+\int_0^t\int_{|\xi|>1}\xi\,N(\d s,\d\xi),\quad t\ge 0
$$ is a $d$-dimensional pure jump L\'evy process
with jump intensity $\nu(\d\xi)$. Given $\sigma$-finite measures $\nu_1$ and $\nu_2$ on $\mbb{R}^d$, we write $\nu_1\land \nu_2 := \nu_1 - (\nu_1-\nu_2)^+ = \nu_2 - (\nu_2-\nu_1)^+$, where the subscript ``$+$'' stands for the upper variation of the signed measure in its Jordan decomposition. A widely used condition in the literature is that for some $x_0 > 0$,
\beqnn
\inf_{x \in \mbb{R}^d,\ |x| \le x_0} [\nu \wedge (\delta_x * \nu) (\mbb{R}^d)]> 0.
\eeqnn
This condition arises from the refined basic coupling method, which
roughly 
requires that the marginal processes can jump together.
However, the above condition is not satisfied when the marginal processes of $(L_t)_{t\ge0}$ are independent, which is exactly the case in our situation.

\item[{\rm(b)}] Studies on the ergodicity of \eqref{OU-proc} typically impose the so-called dissipative conditions. Arapostathis et al. \cite[Corollary 5.2]{APS19}, \cite[Theorem 1.3]{APS22} and Liang et al. \cite[Theorem 1,1]{LMW21} assume
$$
\langle x,b(x)\rangle\le-\lambda|x|^2
$$
for some constant $\lambda>0$ and all large enough $|x|$. Alternatively, Luo and Wang \cite[Theorem 1.1]{LW19} and Majka \cite[Theorem 1.1]{Maj17}, \cite[Corollary 2.7]{Maj19} assume that
$$
\langle x-y,b(x)-b(y)\rangle\le-\lambda|x-y|^2
$$
for some constant $\lambda$ and all $x,y$ with $|x-y|$ large enough. Similar conditions also appear in the study of SDEs driven by Brownian motion, as shown by Eberle \cite[Corollary 2.1]{Eb11}, \cite[Corollary 2]{Eb16}, and in the setting with multiplicative noise, as demonstrated by Zhang and collaborators \cite[Theorems 2.9 and 2.12]{XZ20} and \cite[Theorem 1.2]{ZZ23}.
 However, under our model assumptions, it can be verified that neither of these conditions is satisfied.
To see this, set $z=(x,y)$ with $x,y\ge0$, and consider the drift coefficient
$$
b(z)=(-b_1x^{\alpha_1}+a_1x+\gamma_1, kxy-b_2y^{\alpha_2}+a_2y+\gamma_2)
$$
as in the system \eqref{main sde}. A direct computation gives
$$
\langle z, b(z)\rangle
=-a_1x^{\alpha_1+1}-a_2 y^{\alpha_2+1}
+b_1x^2+kxy^2+b_2y^2+\gamma_1x+\gamma_2y.
$$
The interaction term $kxy^2$ (with $k>0$) 
creates problems even when $\alpha_1,\alpha_2>1$. For example, take $\alpha_1=\alpha_2=1.5$ and set $x=cy$ with some $c>0$. Then as $|z|\rightarrow\infty$, $
\langle z, b(z)\rangle \rightarrow\infty
$ when $k>0$. Hence, $\langle z, b(z)\rangle$ can not be bounded above by a negative quadratic function of $|z|$, 
when $|z|$ is large enough. An analogous issue arises when considering $\langle z-z', b(z)- b(z')\rangle$.

\end{itemize}
\end{rem1}

We next provide a brief introduction to the classical coupling method, and then highlight the advantages of the localized coupling approach. Let $(X_t, Y_t)_{t \ge 0}$ denote the strong solution to \eqref{main sde} with transition semigroup $(P_t)_{t \ge 0}$. We say that a strong Markov process $(X_t, \tilde{X}_t, Y_t, \tilde{Y}_t)_{t \ge 0}$ on $\mathbb{R}_+^4$ is a \textit{Markovian coupling} of the system $(X_t, Y_t)_{t \ge 0}$,  if the marginal process $(\tilde{X}_t, \tilde{Y}_t)_{t \ge 0}$ shares the same transition semigroup $(P_t)_{t \ge 0}$ as $(X_t, Y_t)_{t \ge 0}$.  Denote by $\mbb{P}_{x}(\cdot)$ (resp. $\mbb{P}_{(x, \tilde{x})}(\cdot)$ or $\mbb{P}_{(x,\tilde{x},y,\tilde{y})}(\cdot)$) the law of the precess $(X_t)_{t\ge0}$ (resp. $(X_t,\tilde{X}_t)_{t\ge0}$ or $(X_t,\tilde{X}_t,Y_t,\tilde{Y}_t)_{t\ge0}$) starting from $x$ (resp. $(x,\tilde{x})$ or $(x,\tilde{x}, y, \tilde{y})$). Furthermore, we use $\mbb{E}_x(\cdot)$ (resp. $\mbb{E}_{(x,\tilde{x})}(\cdot)$ or $\mbb{E}_{(x,\tilde{x},y,\tilde{y})}(\cdot)$) to denote the corresponding expectation. The Markovian coupling process \((X_t,\tilde{X}_t, Y_t, \tilde{Y}_t)_{t \ge 0}\) satisfies $$(X_{\tilde{T}+t}, Y_{\tilde{T}+t}) = (\tilde{X}_{\tilde{T}+t}, \tilde{Y}_{\tilde{T}+t})$$ for all $t \ge 0$, where  $$\tilde{T}=\inf\{t \ge 0: (X_t,Y_t)=(\tilde{X}_t,\tilde{Y}_t)\}$$  is called the \textit{succeeding} or \textit{coupling time}. It is essential to minimize the coupling time, since the total variation distance between the marginal distributions of the coupling process \((X_t, \tilde{X}_t, Y_t, \tilde{Y}_t)_{t \ge 0}\) is controlled by the coupling time $\tilde T$ as follows:
\beqnn
\|\delta_{(x,y)} P_t-\delta_{(\tilde{x},\tilde{y})} P_t\|_{\rm Var} \le 2 \mbb{E}_{(x, \tilde{x}, y, \tilde{y})}[\I_{\{(X_t,Y_t)\neq (\tilde{X}_t,\tilde{Y}_t)\}}]=2\mbb{P}_{(x,\tilde{x}, y, \tilde{y})}(\tilde{T}>t)
\eeqnn
for any $x,\tilde{x}, y, \tilde{y}\in\mbb{R}_+.$
Therefore, to obtain \eqref{eq1.10}, it is sufficient to prove that there exist constants $C, \lambda^* > 0$ such that
\beqlb\label{ineq:1.14}
\sup_{x,\tilde{x}, y, \tilde{y}} \mbb{P}_{(x,\tilde{x}, y, \tilde{y})}(\tilde{T}>t) \le C\e^{-\lambda^*t}.
\eeqlb
This classical coupling method is valid for several types of Markov processes; see e.g.,  Chen and Li \cite{CL89}, Eberle \cite{Eb16} and F.Y. Wang \cite{W11}.

 However,
 it is challenging to verify \eqref{ineq:1.14} in our setting, since the associated drift term does not satisfy the so-called dissipative conditions, as noted in Remark \ref{rem for drift term}(b). Roughly speaking, when $X_t$ or $\tilde{X}_t$ becomes large, the predation terms $kX_tY_t$ or $k\tilde{X}_t\tilde{Y}_t$ introduce a strong drift that depends on the current state, which can pull $Y_t$ and $\tilde{Y}_t$ apart and prevent them from meeting. To circumvent this difficulty, we adopt a localized coupling approach. The key idea is to exploit the fact that the Markovian coupling process $(X_t,\tilde{X}_t,Y_t,\tilde{Y}_t)_{t\ge0}$ behaves well when the coupling of the first component $(X_t,\tilde X_t)_{t\ge0}$ is confined to a compact set. Intuitively, instead of aiming for a uniform estimate of the coupling time $\tilde T$, we restrict the coupling process $(X_t,\tilde{X}_t,Y_t,\tilde{Y}_t)_{t\ge0}$ to a smaller subset of $\mathbb{R}^4_+$, and consider the coupling time $\tilde T$ before the coupling process $(X_t,\tilde{X}_t,Y_t,\tilde{Y}_t)_{t\ge0}$ exits this subset. To be more precise, we first let $X_t$ and $\tilde{X}_t$ meet, and estimate 
 the coupling time of  $(X_t)_{t\ge0}$ and $(\tilde X_t)_{t\ge0}$. Then, we restrict $X_t$ and $\tilde{X}_t$ in a bounded set, and consider the coupling time $\tilde T$ before either $X_t$ or $\tilde{X}_t$ exits
 a relatively larger set.
To this end, for any $b\in\mbb{R}_+$, we let $$\tau^+_b=\inf\{t>0:X_t\ge b\},\quad \tau^-_b=\inf\{t>0: X_t\le b\}.$$
Let
$$
\tilde{T}_X=\inf\{t\ge0: X_t=\tilde{X}_t\}
$$
be the coupling time of $(X_t)_{t\ge0}$ and $(\tilde{X}_t)_{t\ge0}$. For all $t\ge0$, we have $X_{t+\tilde{T}_X}=\tilde{X}_{t+\tilde{T}_X}$. It is clear that $\tilde{T} \ge \tilde{T}_X$. The proof of Theorem \ref{main result} is based on the following key observation.

\begin{thm1}\label{main result1}
Suppose that there are some $M,t_0>0$ such that
\beqlb
&\inf\limits_x\mbb{P}_{x}(\tau^-_M<t_0)>0;\label{cond(a)}\\
&\inf\limits_{x \in [0, M], y, \tilde{y}}\mbb{P}_{(x, x, y, \tilde{y})}(\tilde{T}<\tau^+_{2M}\wedge t_0)>0;
\label{cond(b)}\\
&\inf\limits_{x, \tilde{x}}\mbb{P}_{(x,\tilde{x})}(\tilde{T}_X<t_0)>0.\label{cond(c)}
\eeqlb
Then the CBIPC-process $(X_t,Y_t)_{t\ge0}$ to \eqref{main sde} is uniformly ergodic in the total variation distance.
\end{thm1}

If \eqref{cond(a)}, \eqref{cond(b)} and \eqref{cond(c)} are all satisfied, then the conclusion of Theorem \ref{main result} follows directly from Theorem \ref{main result1}. In other words, it suffices to prove the following three propositions.

\bgproposition\label{prop:1.6}
Suppose that $\alpha_1>1$. Then
there exists $M_0>0$ such that \eqref{cond(a)} holds for all $t_0>2$ and $M\ge M_0$.
\edproposition

\bgproposition\label{prop:1.7}
Suppose that $\alpha_i > 1$ for each $i=1, 2$, and either $\sigma_2>0$ or there are constants $\theta_2\in(1,2)$, $C_2>0$ such that
\beqnn
\int_0^x \xi^2\,n_2(\d \xi)\ge C_2x^{2-\theta_2},\quad x\in (0,1].
\eeqnn
Then there exists a constant $M_0>0$ such that for all $M\ge M_0$, 
 there is $t_0:=t_0(M)>0$ so that \eqref{cond(b)} holds.
\edproposition

\bgproposition\label{prop:1.8}
Suppose that $\alpha_1>1$, and either $\sigma_1>0$ or there are constants $\theta_1\in(1,2)$, $C_1>0$ such that
\beqnn
\int_0^x \xi^2\,n_1(\d \xi)\ge C_1x^{2-\theta_1},\quad x\in (0,1].
\eeqnn
Then there is $t_0>0$ so that \eqref{cond(c)} holds.
\edproposition

Indeed, according to Propositions \ref{prop:1.6}, \ref{prop:1.7} and \ref{prop:1.8}, we can see that \eqref{cond(a)}, \eqref{cond(b)} and \eqref{cond(c)} hold simultaneously by taking $t_0$ and $M_0$ to be the maximum of the respective values from these propositions.

\begin{rem1}\label{remark on k le 0}\rm
We emphasize that when $k\le0$, the uniform ergodicity of the process $(X_t,Y_t)_{t\ge0}$ can be proved easily under Condition \ref{Condition}. The details will be provided at the end of Section \ref{sect:3}. Therefore, in the remainder of the article, we primarily discuss the case where $k > 0$.
\end{rem1}

The rest of the paper is organized as follows. In Section \ref{sect:2}, we give the construction of the Markovian coupling process. The proofs of Theorems \ref{main result} and \ref{main result1} are given in Section \ref{sect:3}. The existence and the uniqueness of the strong solution to the stochastic system \eqref{main sde} are presented in the Appendix of this paper.

\section{Markovian coupling for the CBIPC-process}\label{sect:2}

 \setcounter{equation}{0}
 We begin by introducing some convenient notation. For any integer $d\ge1$, let $C^2(\mathbb{R}_+^d)$ be the linear space of twice continuously differentiable functions on $\mathbb{R}_+^d$. Given a function $f\in C^2(\mathbb{R}_+^d)$,
we write
\[
D_{z}f(x)=f(x+z)-f(x)-z\cdot \nabla f(x)
\]
for $x,z\in\mbb{R}_+^d$, where $ \nabla f$ denotes the gradient of $f$ and $\cdot$ is the dot product in $\mbb{R}_+^d$. Now, for any $f\in C^2(\mathbb{R}_+^2)$, we write
\beqlb\label{generator}
Lf(x,y)=L_Xf(x,y)+L_{x,Y}f(x,y),
\eeqlb
where
\beqnn
L_Xf(x,y)=(-b_1x^{\alpha_1}+a_1x+\gamma_1) f'_x(x,y)+\sigma_1xf''_{xx}(x,y)+x\int_0^\infty D_{(\xi,0)}f(x,y)\,n_1(\d \xi)
\eeqnn
and
\beqnn
L_{x,Y}f(x,y)=(kxy-b_2y^{\alpha_2}+a_2y+\gamma_2)f'_y(x,y)
+\sigma_2yf''_{yy}(x,y)
+y\int_0^\infty D_{(0,\xi)}f(x,y)\,n_2(\d \xi),
\eeqnn 
where $f_x'$ denotes the first derivative with respect to $x$, $f_{xx}''$ denotes the second derivative with respect to $x$, and the same for $f'_y$ and $f''_{yy}$. 
Let ${\cal{D}}(L)$ denote the linear space consisting of functions $f\in C^2(\mathbb{R}_+^2)$ such that the integrals on the right-hand side
of \eqref{generator} are convergent and define continuous functions on $\mathbb{R}_+^2$. We shall see that $(L,{\cal{D}}(L))$ is a restriction of the generator of  $(X_t,Y_t)_{t\ge0}$.

\bgproposition\label{prop:2.1}
Let $f\in {\cal{D}}(L)$. Let $A\subset \mbb{R}^2_+$ be such that
\beqnn
\sup_{(x, y)\in A}[|f(x, y)|+|Lf(x, y)|]<\infty,
\eeqnn
and let $
\tau_A=\inf\{t\ge0: (X_t, Y_t)\notin A\}$.
Then for any $t\ge 0$,
\beqlb\label{martingale prob}
f(X_{t\wedge\tau_A},Y_{t\wedge\tau_A})=f(X_0,Y_0)+\int_0^{t\wedge\tau_A}Lf(X_s,Y_s)\,\d s+M_{t\wedge \tau_A},
\eeqlb
where $(M_t)_{t\ge0}$ is a
local
martingale defined by
 \begin{equation}\label{local mart}\begin{split}
M_t =&\sqrt{2\sigma_1}\int_0^{t}\int_0^{X_s}f'_x(X_s,Y_s)\,W_1(\d s,\d u)+
\sqrt{2\sigma_2}\int_0^{t}\int_0^{Y_s}f'_y(X_s,Y_s)\,W_2(\d s, \d u)\cr
&+\int_0^{t}\int_0^\infty\int_0^{X_{s-}}
\left[f(X_{s-}+\xi,Y_{s-})-f(X_{s-},Y_{s-})\right]\,\tilde{N}_1(\d s,\d \xi,\d u)\cr
&+\int_0^{t}\int_0^\infty\int_0^{Y_{s-}}
\left[f(X_{s-},Y_{s-}+\xi)-f(X_{s-},Y_{s-})\right]\,\tilde{N}_2(\d s,\d \xi,\d u).
\end{split}\end{equation}
\edproposition
\proof For any $f\in{\cal{D}}(L)$, we can use It\^{o}'s formula to see that for all $t\ge0$,
\begin{align*}
f(&X_{t\wedge\tau_A},Y_{t\wedge\tau_A})\\
= &
f(X_0,Y_0)
+\int_0^{t\wedge\tau_A}f'_x(X_s,Y_s)
\left[-b_1 X^{\alpha_1}_s + a_1 X_s + \gamma_1\right]\,\d s\\
&
+\int_0^{t\wedge\tau_A}f'_y(X_s,Y_s)\left[k X_sY_s - b_2 Y^{\alpha_2}_s + a_2 Y_s + \gamma_2\right]\,\d s\\
&
+\sqrt{2\sigma_1}\int_0^{t\wedge\tau_A}\int_0^{X_s}f'_x(X_s,Y_s)\,W_1(\d s, \d u)+
\sqrt{2\sigma_2}\int_0^{t\wedge\tau_A}\int_0^{Y_s}f'_y(X_s,Y_s)\,W_2(\d s,\d u)\\
&
+\sigma_1\int_0^{t\wedge\tau_A}f''_{xx}(X_s,Y_s)X_s\,\d s+
\sigma_2\int_0^{t\wedge\tau_A}f''_{yy}(X_s,Y_s)Y_s\,\d s\\
&
+\int_0^{t\wedge\tau_A} \int_0^\infty \int_0^{X_{s-}}
(f(X_s+\xi, Y_s)-f(X_s,Y_s)) \,\tilde{N}_1(\d s, \d \xi, \d u)\\
&
+\int_0^{t\wedge\tau_A} \int_0^\infty \int_0^{Y_{s-}}
(f(X_s, Y_s+\xi)-f(X_s,Y_s)) \,\tilde{N}_2(\d s, \d \xi, \d u)\\
&
+\int_0^{t\wedge\tau_A} \int_0^\infty \int_0^{X_{s-}}
D_{(\xi,0)}f(X_{s-},Y_{s-})\,
 \d s\,n_1(\d \xi)\,\d u\\
&
+\int_0^{t\wedge\tau_A} \int_0^\infty \int_0^{Y_{s-}}
D_{(0,\xi)}f(X_{s-},Y_{s-})\,
\d s\,n_2(\d \xi)\,\d u\\
=&
f(X_0,Y_0)
+\int_0^{t\wedge\tau_A}f'_x(X_s,Y_s)
\left[-b_1 X^{\alpha_1}_s + a_1 X_s + \gamma_1\right]\,\d s\\
&
+\int_0^{t\wedge\tau_A}f'_y(X_s,Y_s)\left[k X_sY_s - b_2 Y^{\alpha_2}_s + a_2 Y_s + \gamma_2\right]\,\d s\\
&
+\sigma_1\int_0^{t\wedge\tau_A}f''_{xx}(X_s,Y_s)X_s\,\d s+
\sigma_2\int_0^{t\wedge\tau_A}f''_{yy}(X_s,Y_s)Y_s\,\d s\\
&
+\int_0^{t\wedge\tau_A}\int_0^\infty X_{s-}
D_{(\xi,0)}f(X_{s-},Y_{s-})\,\d s\,n_1(\d \xi)\\
&
+\int_0^{t\wedge\tau_A}\int_0^\infty Y_{s-}
D_{(0,\xi)}f(X_{s-},Y_{s-})\,\d s\,n_2(\d \xi)\\
&
+\sqrt{2\sigma_1}\int_0^{t\wedge\tau_A}\int_0^{X_s}f'_x(X_s,Y_s)\,W_1(\d s,\d u)+
\sqrt{2\sigma_2}\int_0^{t\wedge\tau_A}\int_0^{Y_s}f'_y(X_s,Y_s)\,W_2(\d s,\d u)\\
&
+\int_0^{t\wedge\tau_A} \int_0^\infty \int_0^{X_{s-}}
(f(X_s+\xi, Y_s)-f(X_s,Y_s))\,
\tilde{N}_1(\d s, \d \xi, \d u)\\
&
+\int_0^{t\wedge\tau_A} \int_0^\infty \int_0^{Y_{s-}}
(f(X_s, Y_s+\xi)-f(X_s,Y_s))\,
\tilde{N}_2(\d s, \d \xi, \d u).
\end{align*}
Then \eqref{martingale prob} holds with $(M_t)_{t\ge0}$ defined by \eqref{local mart}. Under the assumption on $A$ and the definition of $\tau_A$, the integrands in \eqref{local mart} are bounded, which implies that $(M_t)_{t\ge0}$ is a local martingale.\qed

Let $(X_t,Y_t)_{t\ge0}$ and $(\tilde{X}_t, \tilde{Y}_t)_{t\ge0}$ be two nonnegative strong solutions to \eqref{main sde}. Clearly, $(X_t,\tilde{X}_t, Y_t, \tilde{Y}_t)_{t\ge0}$  is a Markovian coupling of the CBIPC-process. We then give a characterization for the generator of the Markovian coupling process in terms of a martingale problem.
Given a function $\Phi\in C^2(\mbb{R}_+^4)$, we write
\begin{align}\label{coupling operator}
\tilde{L}\Phi(x, \tilde{x}, y, \tilde{y})=\tilde{L}_X\Phi(x, \tilde{x}, y, \tilde{y})+\tilde{L}_{x,\tilde{x},Y}\Phi(x,\tilde{x},y,\tilde{y}),
\end{align}
where
\beqnn
\tilde{L}_X\Phi(x, \tilde{x}, y, \tilde{y})
\ar=\ar(-b_1x^{\alpha_1}+a_1x+\gamma_1)\Phi'_x(x, \tilde{x}, y, \tilde{y})
+(-b_1\tilde{x}^{\alpha_1}+a_1\tilde{x}+\gamma_1)\Phi'_{\tilde{x}}(x, \tilde{x}, y, \tilde{y})\cr
\ar\ar
+\sigma_1x\Phi''_{xx}(x, \tilde{x}, y, \tilde{y})
+\sigma_1\tilde{x}\Phi''_{\tilde{x}\tilde{x}}(x, \tilde{x}, y, \tilde{y})
+2\sigma_1(x\wedge\tilde{x})\Phi''_{x\tilde{x}}(x, \tilde{x}, y, \tilde{y})\cr
\ar\ar
+(x\wedge\tilde{x})\int_0^\infty
D_{(\xi,\xi,0,0)}\Phi(x,\tilde{x},y,\tilde{y})\,n_1(\d \xi)\cr
\ar\ar
 +(x-\tilde{x})^+\int_0^\infty D_{(\xi,0,0,0)}\Phi(x,\tilde{x},y,\tilde{y})\,n_1(\d \xi)\cr
 \ar\ar
  +(x-\tilde{x})^-\int_0^\infty D_{(0,\xi,0,0)}\Phi(x,\tilde{x},y,\tilde{y})\,n_1(\d \xi),
 \eeqnn
and
\beqnn
\tilde{L}_{x,\tilde{x},Y}\Phi(x,\tilde{x},y,\tilde{y})
 \ar=\ar(kxy-b_2y^{\alpha_2}+a_2y+\gamma_2)\Phi'_y(x, \tilde{x}, y, \tilde{y})\cr
 \ar\ar+(k\tilde{x}\tilde{y}-b_2\tilde{y}^{\alpha_2}+a_2\tilde{y}+\gamma_2)\Phi'_{\tilde{y}}(x, \tilde{x}, y, \tilde{y})\cr
\ar\ar +\sigma_2y\Phi''_{yy}(x, \tilde{x}, y, \tilde{y})+\sigma_2\tilde{y}\Phi''_{\tilde{y}\tilde{y}}(x, \tilde{x}, y, \tilde{y})+2\sigma_2(y\wedge\tilde{y})\Phi''_{y\tilde{y}}(x, \tilde{x}, y, \tilde{y})\cr
\ar\ar
 +(y\wedge\tilde{y})\int_0^\infty
D_{(0,0,\xi,\xi)}\Phi(x,\tilde{x},y,\tilde{y})\,n_2(\d \xi)\cr
 \ar\ar
 +(y-\tilde{y})^+\int_0^\infty D_{(0,0,\xi,0)}\Phi(x,\tilde{x},y,\tilde{y})\,n_2(\d \xi)\cr
 \ar\ar
 +(y-\tilde{y})^-\int_0^\infty D_{(0,0,0,\xi)}\Phi(x,\tilde{x},y,\tilde{y})\,n_2(\d \xi).
 \eeqnn 
 Here $\Phi''_{x\tilde{x}}$ means the second derivative of the function $\Phi$ with respect to $x$ and $\tilde x$, and the same for $\Phi''_{y\tilde{y}}.$
Let ${\mathcal{D}}(\tilde{L})$ denote the linear space consisting of 
functions $\Phi\in C^2(\mbb{R}_+^4)$ such that the integrals in  \eqref{coupling operator} are convergent and define continuous  functions on $C^2(\mbb{R}_+^4)$. We call $(\tilde{L}, {\mathcal{D}}(\tilde{L}))$ the {\it coupling generator} of the CBIPC-process. The precise meaning of this terminology is made
clear by the martingale problem given in the following.

\bgproposition
Let $\Phi\in {\cal{D}}(\tilde{L})$. Let $A\subset \mbb{R}^4_+$ be such that
\beqnn
\sup_{(x,\tilde{x}, y, \tilde{y})\in A}[|\Phi(x,\tilde{x}, y, \tilde{y})|+|\tilde{L}\Phi(x,\tilde{x}, y, \tilde{y})|]<\infty,
\eeqnn
and let $\tilde \tau_A=\inf\{t\ge0: (X_t, \tilde{X}_t, Y_t, \tilde{Y}_t)\notin A\}$.
Then for any $t\ge 0$,
\beqnn
\Phi(X_{t\wedge\tilde \tau_A}, \tilde{X}_{t\wedge\tilde \tau_A}, Y_{t\wedge\tilde \tau_A},\tilde{Y}_{t\wedge\tilde \tau_A}) - \Phi(X_0,\tilde{X}_0, Y_0,\tilde{Y}_0)- \int_0^{t\wedge\tilde \tau_A}\tilde{L}\Phi(X_s,\tilde{X}_s, Y_s,\tilde{Y}_s)\,\d s
\eeqnn
is a martingale.
\edproposition

\proof The proof is similar to that of Proposition \ref{prop:2.1} and we omit here. \qed

\section{Proof of Theorem \ref{main result}}\label{sect:3}
 \setcounter{equation}{0}

The proof of Theorem \ref{main result} relies essentially on Theorem \ref{main result1}.  We therefore begin by proving the latter.

\noindent \textit{Proof of Theorem $\ref{main result1}$.}\,\,
We first note that $X_t=\tilde X_t$ for all $t\ge 0$ when $X_0=\tilde X_0$.
Then we apply the strong Markov property at $\tau^-_M$, together with  \eqref{cond(a)} and \eqref{cond(b)}, to obtain
\begin{align*}
&\inf_{x,y,\tilde{y}}\mbb{P}_{(x,x,y,\tilde{y})}(\tilde{T}<2t_0)\\
&\ge \inf_{x,y,\tilde{y}}\mbb{P}_{(x,x,y,\tilde{y})}(\tau^-_M<t_0, \tilde{T}\circ\theta_{\tau^-_M}< \tau^+_{2M}\circ\theta_{\tau^-_M} \wedge t_0)\\
&= \inf_{x,y,\tilde{y}}\mbb{E}_{(x, x, y,\tilde{y})}\Big[\I_{\{\tau^-_M<t_0\}}\mbb{E}_{(X_{\tau^-_M}, X_{\tau^-_M}, Y_{\tau^-_M},\tilde{Y}_{\tau^-_M})}(\I_{\{\tilde{T}<\tau^+_{2M}\wedge t_0\}})\Big]\\
&\ge\inf_x\mbb{P}_{x}(\tau^-_M<t_0) \inf_{x\in[0,M], y, \tilde{y}}\mbb{P}_{(x, x, y, \tilde{y})}(\tilde{T}<\tau^+_{2M}\wedge t_0)
>0.\end{align*}
Here and in what follows, $(\theta_t)_{t \ge 0}$ is a collection of  shift operators.  This inequality, together with  the strong Markov property applied at time $\tilde{T}_X$ and \eqref{cond(c)},  gives us
\beqnn
\inf_{x,\tilde{x},y,\tilde{y}}\mbb{P}_{(x,\tilde{x},y,\tilde{y})}(\tilde{T}<3t_0)\ar\ge\ar
\inf_{x,\tilde{x},y,\tilde{y}}\mbb{P}_{(x,\tilde{x},y,\tilde{y})}(T_X<t_0, \tilde{T}\circ\theta_{\tilde{T}_X}<2t_0)\cr
\ar\ge\ar \inf_{x,\tilde{x}}\mbb{P}_{(x,\tilde{x})}(\tilde{T}_X<t_0)\inf_{x,y,\tilde{y}}
\mbb{P}_{(x,x,y,\tilde{y})}(\tilde{T}<2t_0)> 0,
\eeqnn
where in the second inequality we used the fact that $X_{\tilde{T}_X}=\tilde{X}_{\tilde{T}_X}$. Consequently, there exists a constant $\epsilon\in(0,1)$ such that
\beqnn
\sup_{x,\tilde{x},y,\tilde{y}} \mbb{P}_{(x,\tilde{x},y,\tilde{y})}(\tilde{T}\ge 3t_0) \le \epsilon.
\eeqnn

Now applying the Markov property at $3t_0$, we can establish that
\beqnn
\sup_{x,\tilde{x},y,\tilde{y}}\mbb{P}_{(x,\tilde{x},y,\tilde{y})}(\tilde{T}\ge 6t_0)
\ar=\ar\sup_{x,\tilde{x},y,\tilde{y}}\mbb{E}_{(x,\tilde{x},y,\tilde{y})}\left[\mbb{P}_{(x,\tilde{x},y,\tilde{y})}(\tilde{T}\ge 6t_0
|{\cal{F}}_{3t_0})\right]\cr
\ar=\ar\sup_{x,\tilde{x},y,\tilde{y}}\mbb{E}_{(x,\tilde{x},y,\tilde{y})}\left[\I_{\{\tilde{T}\ge 3t_0\}}\mbb{P}_{(X_{3t_0},\tilde{X}_{3t_0},Y_{3t_0},\tilde{Y}_{3t_0})}(\tilde{T}\ge 3t_0)\right]\cr
\ar\le\ar\sup_{x,\tilde{x},y,\tilde{y}}\mbb{E}_{(x,\tilde{x},y,\tilde{y})}\left[\I_{\{\tilde{T}\ge 3t_0\}}\sup\limits_{x,\tilde{x},y,\tilde{y}}\mbb{P}_{(x,\tilde{x},y,\tilde{y})}(\tilde{T}\ge 3t_0)\right]\le \epsilon^2.
\eeqnn
Hence, by induction,
$$
\sup\limits_{x,\tilde{x},y,\tilde{y}}\mbb{P}_{(x,\tilde{x},y,\tilde{y})}(\tilde{T}\ge 3nt_0)\le \epsilon^n,\qquad n=1,2,\dots.
$$

Note that for any fixed $t > 3t_0$, there exists an integer $n_0 \ge 1$ such that $$3 n_0 t_0 < t \le 3(n_0+1) t_0.$$ Then, for such $t$ we have
\begin{align*}
\sup\limits_{x,\tilde{x},y,\tilde{y}}\mbb{P}_{(x,\tilde{x},y,\tilde{y})}(\tilde{T}\ge t) \le
\sup\limits_{x,\tilde{x},y,\tilde{y}}\mbb{P}_{(x,\tilde{x},y,\tilde{y})}(\tilde{T}\ge 3n_0t_0)
\le
\epsilon^{n_0}
\le
\mrm{e}^{-\lambda_*t},
\end{align*}
where $\lambda_*:=-(3t_0)^{-1}\ln\epsilon$.
On the other hand, since the probability is clearly bounded by $1$ when $t \le 3t_0$, we conclude that for all $t \ge 0$,
\beqnn
\sup\limits_{x,\tilde{x},y,\tilde{y}}\mbb{P}_{(x,\tilde{x},y,\tilde{y})}(\tilde{T}\ge t)\le \left(1 \vee \mrm{e}^{3\lambda_*t_0}\right) \mrm{e}^{-\lambda_*t}.
\eeqnn
Letting $C = 2 \left(1 \vee \e^{3\lambda_* t_0}\right)$, we obtain the following uniform exponential contraction:
\beqlb\label{ineq:3.2}
\begin{split}
\|\delta_{(x,y)} P_t-\delta_{(\tilde{x},\tilde{y})} P_t\|_{\rm Var} &\le 2 \I_{\{(x,y)\neq (\tilde{x},\tilde{y})\}}\left[\sup_{x,\tilde{x},y,\tilde{y}}\mbb{P}_{(x,\tilde{x},y,\tilde{y})}(\tilde{T} \ge t)\right]\\
&
\le
C \mrm{e}^{-\lambda_*t} \I_{\{(x,y)\neq (\tilde{x},\tilde{y})\}},\quad t \ge 0, x,\tilde x, y, \tilde y \ge0.\end{split}
\eeqlb

Having established the uniform exponential contraction in \eqref{ineq:3.2}, we now proceed to prove the uniform ergodicity.  By the convexity of the total variation distance, we have for any $\gamma, \eta \in \mcr{P}(\mbb{R}^2_+)$ and $\pi \in \mcr{C}(\gamma,\eta)$,
\beqnn
\|\gamma P_t-\eta P_t\|_{\rm Var}
\le
\int_{\mbb{R}_+^2\times\mbb{R}_+^2}
\|\delta_{z}P_t-\delta_{\tilde{z}}P_t\|_{\rm Var}\, \pi(\d z,\d \tilde{z}),
\eeqnn  where $\mcr{C}(\gamma,\eta)$ denotes the set of all coupling measures on $\mbb{R}^4_+$ of $\gamma$ and $\eta$;
see, e.g., Villani \cite[Theorem 4.8]{V09}. Hence, it follows from \eqref{ineq:3.2} that
\beqlb\label{ineq:3.3}
\|\gamma P_t-\eta P_t\|_{\rm Var}\le C \mrm{e}^{-\lambda_*t}\|\gamma-\eta\|_{{\rm Var}},\quad t\ge0,\,\, \gamma,\eta\in  \mcr{P}(\mbb{R}^2_+).
\eeqlb
Then, for sufficiently large $r>0$, the operator $P_r$
is contractive on $(\mcr{P}(\mbb{R}^2_+), \|\cdot\|_{\rm Var})$. By Cohn \cite[Proposition 4.1.8]{Co13}, $(\mcr{P}(\mbb{R}^2_+), \|\cdot\|_{\rm Var})$ is a complete metric
space. Then, according to 
the Banach fixed point theorem, there is a unique $\gamma_r\in\mcr{P}(\mbb{R}^2_+)$ such that $\gamma_rP_r=\gamma_r$. Fix such an $r>0$ and define $\gamma=r^{-1}\int_0^r\gamma_rP_s\,\d s$. Clearly, $\gamma \in \mcr{P}(\mbb{R}^2_+)$. For any $0\le t <r$, applying the  Chapman-Kolmogorov equation, we see that
\beqlb\label{ineq:3.4}
\begin{split}
\gamma P_t &= r^{-1}\int_0^r\gamma_rP_{t+s}\,\d s= r^{-1}\int_t^{r+t}\gamma_rP_s\,\d s\\
&=
r^{-1}\int_t^r\gamma_rP_s\,\d s
+
r^{-1}\int_0^t\gamma_rP_{r+s}\,\d s=
r^{-1}\int_t^r\gamma_rP_s\,\d s
+
r^{-1}\int_0^t\gamma_rP_s\,\d s\\
&=
r^{-1}\int_0^r\gamma_rP_s\,\d s=\gamma.\end{split}
\eeqlb

Now for any arbitrary $t\ge r$, let $k\in\mbb{N}$ be such that $0\le t-kr<r$. By the semigroup property and the definition of $\gamma$, we have $\gamma P_{kr}=\gamma$. It follows that
$$
\gamma P_t=\gamma P_{kr+(t-kr)}=\gamma P_{kr}P_{t-kr}
=\gamma P_{t-kr}.
$$
Since
$0 \le t-kr <r$, \eqref{ineq:3.4} implies that $\gamma P_{t-kr}=\gamma$, which completes the verification that $$
\gamma P_t=\gamma
$$
for all $t\ge0$. Consequently, it follows from \eqref{ineq:3.3} that
$$
\|\gamma-\eta P_t\|_{\rm Var}\le C \mrm{e}^{-\lambda_* t}\|\gamma-\eta\|_{\rm Var},\quad t \ge 0, \eta\in \mcr{P}(\mbb{R}^2_+),
$$
which establishes the uniform ergodicity of the process $(X_t,Y_t)_{t\ge0}$.\qed

\subsection{Proof of Proposition \ref{prop:1.6}}

\noindent{\it Proof of Proposition $\ref{prop:1.6}$.}\,\, Since $\alpha_1>1$, we can define a bounded nonnegative function
$$
g(x)=2-(1+x)^{-{(\alpha_1-1)}/{2}},\qquad x\ge 0.
$$
It is easy to verify that  for any $x>0$,
\beqnn
g'(x)=(\alpha_1-1)(1+x)^{-{(\alpha_1+1})/{2}}/2>0\quad\mbox{and}\quad g''(x)=(1-\alpha_1^2)(1+x)^{-{(\alpha_1+3)}/{2}}/4<0.
\eeqnn
We choose a sufficiently large $M_0>0$ such that for all $M\ge M_0$ and $x\ge M$,
$$
-b_1x^{\alpha_1}+a_1x+\gamma_1
\le-2(1+x)^{(\alpha_1+1)/2}/(\alpha_1-1).
$$
Substituting these estimates into \eqref{generator} and applying the mean-value theorem, we obtain that for all $x \ge M$,
\beqlb\label{ineq:3.5}
{L}_Xg(x)\le(-b_1x^{\alpha_1}+a_1x+\gamma_1)g'(x)\le -1.
\eeqlb
By Proposition \ref{prop:2.1}, we have
$$
g(X_{t\wedge\tau^-_M})=
g(X_0)+\int_0^{t\wedge\tau^-_M}L_Xg(X_s)\,\d s+M_{t\wedge\tau^-_M}£¬
$$
where $(M_{t\wedge\tau^-_M})_{t\ge0}$ is a martingale. Then, taking expectations on both sides, we obtain that for all $x\ge M$,
\beqnn
0\le \mbb{E}_x\left[g(X_{t\wedge\tau^-_M})\right]=
g(x)+\mbb{E}_{x}\left[\int_0^{t\wedge\tau^-_M} {L}_Xg(X_s)\,\d s\right]\le
g(x)-\mbb{E}_x[t\wedge \tau^-_M],
\eeqnn
where the last inequality is the result of \eqref{ineq:3.5}. Now, letting $t\to\infty$ on both sides, we deduce that
\beqnn
\sup_{x \ge M}\mbb{E}_x[\tau^-_M]\le \sup_{x \ge M}g(x)=2.
\eeqnn
Finally, by applying the Markov inequality,  for any $t_0>2$, we have
$$
\sup_{x \ge M}\mbb{P}_x(\tau^-_M\ge t_0)\le \frac{\sup_{x \ge M}\mbb{E}_x[\tau^-_M]}{t_0}<1,
$$
which implies
$$
\inf_x\mbb{P}_x(\tau^-_M < t_0)= \inf_{x \ge M}\mbb{P}_x(\tau^-_M < t_0)>0.
$$
The proof is complete.
 \qed

\subsection{Proof of Proposition \ref{prop:1.7}}\label{sect:3.2}

Let $(X_t, {X}_t, Y_t,\tilde{Y}_t)_{t \ge 0}$ be the Markovian coupling process  with $X_0 = \tilde{X}_0 \ge 0$ and $Y_0, \tilde{Y}_0\ge0$. In particular, due to the characterization of \eqref{main sde}, we have $X_t=\tilde X_t$ for all $t\ge0$. 
The proof of Proposition \ref{prop:1.7} is more complex than that of Proposition \ref{prop:1.6}. We thus give a brief overview of the argument. First, in Section \ref{Section3.2.1} we establish some local comparison principles. These allow us to control the difference process $(Y_t-\tilde{Y}_t)_{t\ge0}$ by an auxiliary process $(\bar{Z}_t)_{t\ge0}$, which is constructed to be independent of the first component process $(X_t)_{t\ge0}$ until the latter exits a relatively larger set.  Second, in Section \ref{Section3.2.2}, we derive tail estimates for both the hitting time of $(\bar{Z}_t)_{t\ge0}$ to zero and the exit time $\tau^+_{2M}$ of $(X_t)_{t\ge0}$. These estimates are obtained through Lyapunov function arguments and some probabilistic ideas. Finally, in Section \ref{Section3.2.3}, using these results in previous sections together with the strong Markov property and the independence, we complete the proof of Proposition \ref{prop:1.7}.

\subsubsection{{\bf Local comparison principles}}\label{Section3.2.1}  We start by presenting a technical construction which we will use later. For each integer $n \ge 0$, define $a_n=\mathrm{exp}(-n(n+1)/2)$. Let $\psi_n $ be a nonnegative  continuous function supported on $(a_n,a_{n-1})$ such that $$
\int_{a_n}^{a_{n-1}}\psi_n(x)\,\d x=1
$$
and $\psi_n(x)\le2/(nx)$ for every $x>0$. For $n \ge 0$ and $z\in\mathbb{R}$, let
$$
\phi_n(z)=\int_0^{z}\,\d y\int_0^y\psi_n(x)\,\d x.
$$
For any $z\in\mbb{R}$, we write $z^+=z\vee0$ for its positive part.

\begin{lem1}\label{le:3.1}
The sequence of the functions $\{\phi_n(x)\}_{n\ge1}$ satisfies
\begin{itemize}
\item[{\rm(1)}] for all $z\in\mathbb{R}$, $\phi_n(z)\rightarrow z^+$ increasingly as $n \rightarrow \infty;$

\item[{\rm(2)}] for all $z,\xi \le 0$, $\phi_n(z)=0$ and $D_\xi \phi_n(z)=0;$

\item[{\rm(3)}]for all $z>0$, $0<\phi'_n(z) \le 1$ and
$0\le z\phi''_n(z)\le2/n;$

\item[{\rm(4)}]for $z, \xi \ge 0$,
$0 \le z D_\xi \phi_n(z) \le(1+2z)[\xi\wedge(\xi^2/n)].$\end{itemize}
\end{lem1}

\proof For the proof of assertion {\rm (1)}, we define \[
F_n(y)=\int_0^y\psi_n(x)\,\d x,\quad y\in \mbb{R}.\] Using the fact that $\psi_n$ is supported on $(a_n, a_{n-1})$ and $a_n>0$, we have for all $z\le 0$
$$\phi_n(z)=F_n(z)=0.$$
This implies that the limit holds for $z\le0$.
We next consider the case of $z> 0.$ Fix any $y>0$.
 Since $a_n\to 0$ as $n \to \infty$,  we have
$y\ge a_{n-1}$ for all sufficiently large $n$, and hence
$F_n(y)=1$. Moreover, for any $n\ge1$ and $y>0$, we have $F_n(y)\in[0,1]$. Then for each $z>0$, by the dominated convergence theorem,
$$
\phi_n(z)=\int_0^z F_n(y)\,\d y\rightarrow\int_0^z 1\,\d y = z\quad\mbox{as}\quad n\rightarrow\infty.
$$

The proof of (2) is trivial since $\psi_n(x)=0$ for $ x\le 0.$

An direct application of the properties of $\psi_n$ yields (3).

We finally turn to prove (4).
On the one hand, by the fact that $\phi'_n\le 1$ and the mean-value theorem, one sees that
$$
D_\xi \phi_n(z) \le \phi_n(z+\xi)-\phi_n(z) + \xi\phi'_n(z) \le 2\xi.
$$
On the other hand, by Taylor's expansion, for any $z,\xi\ge0$, there is an $\eta\in (z, z+\xi)$ such that
\beqnn
z D_\xi \phi_n(z) = z
\phi''_n(\eta) \xi^2/2
\le \eta \phi''_n(\eta) \xi^2/2 \le \xi^2/n,
\eeqnn
where we used the assertion (3)
in the last inequality.
\qed

Using on the technical construction above, we first  establish a local order-preserving property for $(Y_t,\tilde{Y}_t)_{t\ge0}$.
Recall that
$$
\tau^+_{2M}=\inf\{t>0: X_t\ge 2M\}.
$$
\begin{lem1}\label{le:3.2}
Given any $M>0$, it holds that
$$
\inf_{x\in[0,M], y>\tilde{y}\ge0}\mbb{P}_{(x,x,y,\tilde{y})}\bigl(Y_t \ge \tilde{Y}_t~\text{for~all~} t\in[0,\tau^+_{2M}]\bigr)=1.
$$
\end{lem1}

\proof Assuming that $Y_0 > \tilde{Y}_0 \ge 0$, we define $\Delta_t=\tilde{Y}_t-Y_t$ for every $t\ge0$. For $N>0$, we define the stopping time $\varsigma_N=\inf\{t>0: |\Delta_t|\ge N\}$.  From \eqref{main sde}, we obtain
\beqnn
\Delta_{t\wedge\tau^+_{2M}\wedge\varsigma_N}\ar=\ar \Delta_0+
\int_0^{t\wedge\tau^+_{2M}\wedge \varsigma_N} \left[(k X_s+a_2)\Delta_{s-}- b_2 (\tilde{Y}^{\alpha_2}_s-Y^{\alpha_2}_s)\right]\,\d s\cr
\ar\ar
+\sqrt{2\sigma_2}
\int_0^{t\wedge\tau^+_{2M}\wedge \varsigma_N}\int_{Y_s}^{\tilde{Y}_s}\,W_2(\d s,\d u)\cr
\ar\ar+
\int_0^{t\wedge\tau^+_{2M}\wedge \varsigma_N}\int_0^\infty\int_{Y_{s-}}^{\tilde{Y}_{s-}}\I_{\{\Delta_{s-} > 0\}}\xi\,\tilde{N}_2(\d s,\d \xi,\d u)\cr
\ar\ar
-\int_0^{t\wedge\tau^+_{2M}\wedge \varsigma_N}\int_0^\infty\int_{\tilde{Y}_{s-}}^{Y_{s-}}\I_{\{\Delta_{s-} \le 0\}}\xi\,\tilde{N}_2(\d s,\d \xi, \d u).
\eeqnn
Let $\{\phi_n\}_{n\ge1}$ be the sequence of the functions given in Lemma \ref{le:3.1}. Since $y > \tilde y$, $\phi_n(\Delta_0)=0$ for all $n\ge1$. Applying It\^{o}'s formula together with Lemma \ref{le:3.1}(2) yields
\beqnn
\phi_n(\Delta_{t\wedge\tau^+_{2M}\wedge \varsigma_N})
\ar=\ar \int_0^{t\wedge\tau^+_{2M}\wedge \varsigma_N} \left[(k X_s+a_2)\Delta_{s-}- b_2 (\tilde{Y}^{\alpha_2}_s-Y^{\alpha_2}_s)\right] \phi'_n(\Delta_{s-})\,\d s\cr
\ar\ar+\int_0^{t\wedge\tau^+_{2M}\wedge \varsigma_N} \Delta_{s-}\I_{\{\Delta_{s-}>0\}}\,\d s \int_0^\infty
D_\xi \phi_n(\Delta_{s-})\,n_2(\d \xi)\cr
\ar\ar+\int_0^{t\wedge\tau^+_{2M}\wedge \varsigma_N} (-\Delta_{s-})\I_{\{\Delta_{s-}\le 0\}}\,\d s \int_0^\infty
D_{-\xi} \phi_n(\Delta_{s-})\,n_2(\d \xi)\cr
\ar\ar+\sigma_2\int_0^{t\wedge\tau^+_{2M}\wedge \varsigma_N}\phi''_n(\Delta_{s-})\Delta_{s-}\,\d s+M_t\cr
\ar=\ar \int_0^{t\wedge\tau^+_{2M}\wedge \varsigma_N} \left[(k X_s+a_2)\Delta_{s-}- b_2 (\tilde{Y}^{\alpha_2}_s-Y^{\alpha_2}_s)\right] \phi'_n(\Delta_{s-}) \I_{\{\Delta_{s-}>0\}}\d s\cr
\ar\ar+\int_0^{t\wedge\tau^+_{2M}\wedge \varsigma_N} \Delta_{s-}\I_{\{\Delta_{s-}>0\}}\,\d s \int_0^\infty
D_\xi \phi_n(\Delta_{s-})\,n_2(\d \xi)\cr
\ar\ar
+\sigma_2\int_0^{t\wedge\tau^+_{2M}\wedge \varsigma_N} \phi''_n(\Delta_s)\Delta_{s-}\I_{\{\Delta_{s-}>0\}}\,\d s
+M_t,
\eeqnn where $(M_t)_{t\ge0}$ is a martingale.
Thanks to Lemma \ref{le:3.1}(3) and (4) and the fact that
$$
(\tilde{Y}_s^{\alpha_2}-Y_s^{\alpha_2})\I_{\{\Delta_{s-}\ge0\}}\ge(\tilde{Y}_s-Y_s)^{\alpha_2}\I_{\{\Delta_{s-}\ge0\}}\ge 0,
$$
we additionally have
\beqnn
\phi_n(\Delta_{t\wedge\tau^+_{2M}\wedge \varsigma_N})
\ar\le\ar
\int_0^{t\wedge\tau^+_{2M}\wedge \varsigma_N} (2kM+a^+_2)\Delta_{s-}\I_{\{\Delta_{s-}>0\}}\,\d s\cr
\ar\ar+\int_0^{t\wedge\tau^+_{2M}\wedge \varsigma_N} (1+2\Delta_{s-})\I_{\{\Delta_{s-}>0\}}\,\d s \int_0^\infty
(\xi\wedge \xi^2 n^{-1})\,n_2(\d \xi)\cr
\ar\ar+2\sigma_2n^{-1}\int_0^{t\wedge\tau^+_{2M}\wedge \varsigma_N} \I_{\{\Delta_{s-}>0\}}\,\d s
+M_t.
\eeqnn
Then, taking expectations on both sides of the equation above, we obtain
\beqnn
\mbb{E}_{(x,x,y,\tilde{y})}[\phi_n(\Delta_{t\wedge\varsigma_N\wedge\tau^+_{2M}})]
\ar\le\ar  (2kM+a^+_2)\int_0^t\mbb{E}_{(x,x,y,\tilde{y})}[\Delta^+_{s\wedge\varsigma_N\wedge\tau^+_{2M}}]\,\d s\cr
\ar\ar+\int_0^\infty(\xi\wedge \xi^2 n^{-1})\,n_2(\d \xi)\int_0^t\mbb{E}_{(x,x,y,\tilde{y})}
[1+2\Delta^+_{s\wedge\varsigma_N\wedge\tau^+_{2M}}]\,\d s
\cr
\ar\ar
+2\sigma_2n^{-1}\mbb{E}_{(x,x,y,\tilde{y})}[t\wedge\varsigma_N\wedge\tau^+_{2M}].
\eeqnn
Now, letting $n\rightarrow\infty$, we deduce from
Lemma \ref{le:3.1}(1) that
$$
\mbb{E}_{(x,x,y,\tilde{y})}[\Delta^+_{t\wedge\zeta_N\wedge\tau^+_{2M}}]\le (2kM+a^+_2)\int_0^t\mbb{E}_{(x,x,y,\tilde{y})}
[\Delta^+_{s\wedge\varsigma_N\wedge\tau^+_{2M}}]\,\d s.
$$
The Gronwall inequality implies that for all $t\ge0$, $\mbb{E}_{(x,x,y,\tilde{y})}[\Delta^+_{t\wedge\varsigma_N\wedge\tau^+_{2M}}]=0.$ Finally, letting $N\rightarrow\infty$ and applying Fatou's lemma, we conclude that
$\mbb{E}_{(x,x,y,\tilde{y})}[\Delta^+_{t\wedge\tau^+_{2M}}]=0$ for all $t \ge 0$, which yields the desired result.\qed

By Lemma \ref{le:3.2}, for any $M\ge0$, whenever the first component $(X_t,\tilde X_t)_{t\ge0}$ of the Markovian coupling process coincides and does not exceed $2M$, $(Y_t,\tilde{Y}_t)_{t\ge0}$ is locally order-preserving. Then, in order to consider the second component $(Y_t,\tilde Y_t)$ for all $\tilde T_X\le t\le \tau^+_{2M}$, it is enough to consider only the case that $Y_0>\tilde{Y}_0\ge0$. The case that $\tilde{Y}_0> Y_0\ge0$ can be handled by symmetry.
We next introduce two auxiliary processes. Given the process $(X_t)_{t\ge0}$, we consider the following SDE:
\begin{equation}\label{Z sde}\begin{split}
Z_t=& Y_0-\tilde{Y}_0 + \int_0^t (kX_sZ_s - b_2 Z_s^{\alpha_2} + a_2 Z_s)\, \d s +\sqrt{2\sigma_2}\int_0^t \int_0^{Z_s}\,W_2(\d s,\d u) \\
&
+ \int_0^t \int_0^\infty \int_0^{Z_{s-}} \xi  \,\tilde{N}_2(\d s, \d \xi, \d u),
\end{split}\end{equation}
where the  coefficients and the noises are consistent with these in \eqref{main sde}. By Proposition \ref{unique strong solution}, there exists a unique nonnegative strong solution $(X_t, Z_t)_{t\ge0}$. Define the following hitting time
$$
\zeta_0=\inf\{t\ge0: Z_t=0\}.
$$
By \eqref{Z sde}, $Z_{t+\zeta_0}=0$ for all $t\ge0$. Recall again that
$$
\tilde{T}=\inf\{t>0: (X_t,Y_t)=(\tilde{X}_t,\tilde{Y}_t)\} \quad \text{and}  \quad \tau^+_{2M}=\inf\{t>0: X_t\ge 2M\}.
$$

In view of Lemma \ref{le:3.2}, it is meaningful to consider
the following two local comparison principles.

\begin{lem1}\label{le:3.3}
 For any $M>0$, it holds that
 $$
\inf_{x\in [0,M],y>\tilde{y}\ge0}\mbb{P}_{(x,x,y,\tilde{y})}
\left(\,
Y_t-\tilde{Y}_t\le Z_t~~\text{for~all}~~ t\in[0,\tau^+_{2M}]\,
\right)=1.
 $$
In particular, $$\inf_{x\in [0,M],y>\tilde{y}\ge0}\mbb{P}_{(x,x,y,\tilde{y})}(\tilde{T}\wedge\tau^+_{2M} \le \zeta_0\wedge\tau^+_{2M})=1.$$
\end{lem1}

\proof Define $\eta_t=Y_t-\tilde{Y}_t-Z_t$ for $t\ge0$. For $N>0$, we define the stopping time $\iota_N=\inf\{t>0: |\eta_t|\ge N\}$. From \eqref{main sde} and \eqref{Z sde}, we obtain
\beqnn
\eta_{t\wedge\tau^+_{2M}\wedge\iota_N}\ar=\ar
\int_0^{t\wedge\tau^+_{2M}\wedge\iota_N} \left[(k X_s+a_2)\eta_s- b_2 (Y^{\alpha_2}_s-\tilde{Y}^{\alpha_2}_s-Z_s^{\alpha_2})\right]\,\d s\cr
\ar\ar
+\sqrt{2\sigma_2}
\int_0^{t\wedge\tau^+_{2M}\wedge\iota_N}\int_{Z_s}^{Y_s-\tilde{Y}_s}\I_{\{\eta_s>0\}}\,W_2(\d s,\d u)\cr
\ar\ar+
\int_0^{t\wedge\tau^+_{2M}\wedge\iota_N}\int_0^\infty\int_{Z_{s-}}^{Y_{s-}-\tilde{Y}_{s-}}\I_{\{\eta_{s-}>0\}}\xi\,\tilde{N}_2(\d s,\d \xi,\d u)\cr
\ar\ar
-
\int_0^{t\wedge\tau^+_{2M}\wedge\iota_N}\int_0^\infty\int_{Y_{s-}-\tilde{Y}_{s-}}^{Z_{s-}}\I_{\{\eta_{s-}\le0\}}\xi\,\tilde{N}_2(\d s,\d \xi, \d u).
\eeqnn
Let $\{\phi_n\}$ be the sequence of the functions given in Lemma \ref{le:3.1}. Applying It\^{o}'s formula together with Lemma \ref{le:3.1}(2) yields
\beqnn
\phi_n(\eta_{t\wedge\tau^+_{2M}\wedge \iota_N})
\ar=\ar \int_0^{t\wedge\tau^+_{2M}\wedge \iota_N} \left[(k X_s+a_2)\eta_{s-}- b_2 (Y^{\alpha_2}_s-\tilde{Y}^{\alpha_2}_s-Z^{\alpha_2}_s)\right] \phi'_n(\eta_{s-})\,\d s\cr
\ar\ar+\int_0^{t\wedge\tau^+_{2M}\wedge \iota_N} \eta_{s-}\I_{\{\eta_{s-}>0\}}\,\d s \int_0^\infty
D_\xi \phi_n(\eta_{s-})\,n_2(\d \xi)\cr
\ar\ar+\int_0^{t\wedge\tau^+_{2M}\wedge \iota_N} (-\eta_{s-})\I_{\{\eta_{s-}\le 0\}}\,\d s \int_0^\infty
D_{-\xi} \phi_n(\eta_{s-})\,n_2(\d \xi)\cr
\ar\ar+\sigma_2\int_0^{t\wedge\tau^+_{2M}\wedge \iota_N}\phi''_n(\eta_{s-})\eta_{s-}\,\d s+ M_t\cr
\ar=\ar \int_0^{t\wedge\tau^+_{2M}\wedge \iota_N} \left[(k X_s+a_2)\eta_{s-}- b_2 (Y^{\alpha_2}_s-\tilde{Y}^{\alpha_2}_s-Z^{\alpha_2}_s)\right] \phi'_n(\eta_{s-}) \I_{\{\eta_{s-}>0\}}\,\d s\cr
\ar\ar+\int_0^{t\wedge\tau^+_{2M}\wedge \iota_N} \eta_{s-}\I_{\{\eta_{s-}>0\}}\,\d s \int_0^\infty
D_\xi \phi_n(\eta_{s-})\,n_2(\d \xi)\cr
\ar\ar
+\sigma_2\int_0^{t\wedge\tau^+_{2M}\wedge \iota_N} \phi''_n(\eta_s)\eta_{s-}\I_{\{\eta_{s-}>0\}}\,\d s
+M_t,
\eeqnn where $(M_t)_{t\ge0}$ is a martingale.
Thanks to Lemma \ref{le:3.1}(3) and (4), and the fact that
$$
(Y_s^{\alpha_2}-\tilde{Y}_s^{\alpha_2}-Z_s^{\alpha_2})\I_{\{\eta_{s-}>0\}}>
(Y_s^{\alpha_2-1}(Y_s-\tilde{Y}_s)-Z_s^{\alpha_2})\I_{\{\eta_{s-}>0\}}
>[(Y_s-\tilde{Y}_s)^{\alpha_2}-Z_s^{\alpha_2}]\I_{\{\eta_{s-}>0\}}>0,
$$
we additionally have
\beqnn
\phi_n(\eta_{t\wedge\tau^+_{2M}\wedge \iota_N})
\ar\le\ar
\int_0^{t\wedge\tau^+_{2M}\wedge \iota_N} (2kM+a^+_2)\eta_{s-}\I_{\{\eta_{s-}>0\}}\,\d s\cr
\ar\ar+\int_0^{t\wedge\tau^+_{2M}\wedge \iota_N} (1+2\eta_{s-})\I_{\{\eta_{s-}>0\}}\,\d s \int_0^\infty
(\xi\wedge \xi^2 n^{-1})\,n_2(\d \xi)\cr
\ar\ar+2\sigma_2n^{-1}\int_0^{t\wedge\tau^+_{2M}\wedge \iota_N} \I_{\{\eta_{s-}>0\}}\,\d s
+M_t.
\eeqnn
The proof of the remaining part is similar to that of Lemma \ref{le:3.2}, so it is omitted here. \qed

Next, for any $M>0$, let us consider another SDE:
\beqlb\label{bar Z sde}\begin{split}
\bar{Z}_t=& Y_0-\tilde{Y}_0 + \int_0^t (2kM\bar{Z}_s - b_2 \bar{Z}_s^{\alpha_2} + a_2 \bar{Z}_s) \,\d s  \\
&
+\sqrt{2\sigma_2}\int_0^t \int_0^{\bar{Z}_s}\,W_2(\d s,\d u) + \int_0^t \int_0^\infty \int_0^{\bar{Z}_{s-}} \xi  \,\tilde{N}_2(\d s, \d \xi, \d u).\end{split}
\eeqlb
By Dawson and Li \cite[Theorem 2.5]{DaL12}, there 
exists a unique nonnegative strong solution $(\bar{Z}_t)_{t\ge0}$.
Define  the following hitting times
$$
\bar{\zeta}_0=\inf\{t\ge0: \bar{Z}_t=0\}, \quad \bar{\zeta}_{1/n}=\inf\{t\ge0: \bar{Z}_t\le 1/n\},\,\, n\ge1.
$$
Clearly, $\lim_{n\rightarrow\infty}\bar{\zeta}_{1/n}=\bar{\zeta}_0$. By \eqref{bar Z sde}, $\bar{Z}_{ t+\bar{\zeta}_0}=0$ for any $t\ge0$.

\begin{lem1}\label{le:3.4}
 For any $M>0$,  it holds that
 $$\inf_{x\in [0,M], y>\tilde{y}\ge0}
\mbb{P}_{(x,x,y,\tilde{y})}
\Big(\,
Z_t\le \bar{Z}_{t}~~\text{for~all}~~ t\in[0,\tau^+_{2M}]\,
\Big)=1.
 $$
 In particular,
 $$\inf_{x\in [0,M], y>\tilde{y}\ge0}\mbb{P}_{(x,x,y,\tilde{y})} ( \zeta_0\wedge\tau^+_{2M} \le \bar{\zeta}_0)=1.$$
\end{lem1}

\proof The proof is quite similar to that of Lemma \ref{le:3.2}, and so it is omitted here. \qed

\subsubsection{{\bf Tail estimates of stopping times $\bar{\zeta}_0$ and $\tau^+_{2M}$}}\label{Section3.2.2} We need to obtain tail estimates for the stopping times $\bar{\zeta}_0$ and $\tau^+_{2M}$. Let us first consider  $\bar{\zeta}_0$. For $z\ge0$ and $f\in C^2(\mbb{R}_+)$,   write
\beqlb\label{generator of bar Z}
L_{\bar{Z}}f(z) \ar=\ar [(2kM+a_2)z-b_2z^{\alpha_2}]f'(z)
+\sigma_2zf''(z)
+z\int_0^\infty D_\xi f(z)\,n_2(\d \xi).
\eeqlb
Let ${\mathcal{D}}(L_{\bar{Z}})$ denote the linear space consisting of functions $f\in C^2(\mbb{R}_+)$ such that the integral on the right-hand side of \eqref{generator of bar Z} is convergent and defines a continuous function on $\mbb{R}_+$. Following the proof of Proposition \ref{prop:2.1}, we can show that $(L_{\bar{Z}}, {\mathcal{D}}(L_{\bar{Z}}))$ is a restriction of the generator of the process $(\bar{Z}_t)_{t\ge0}$ to \eqref{bar Z sde}.

We define a function
under Condition \ref{Condition} by
\beqlb\label{f0}
f(z) =
\begin{cases}
	1 + z^\beta,&\quad z\in[0,l_0],\\
	f(l_0)+\frac{\beta l_0^\beta}{1-\beta}(1-\mrm{e}^{-(1-\beta)(z-l_0)/l_0}), & \quad z\in(l_0,l_1],\\
f(l_1)+c_0\int_0^{z-l_1}\frac{2}{b_2(c_1s+l_1)^{\alpha_2}}\,\d s, &\quad z\in(l_1,\infty),
\end{cases}
\eeqlb
where
\beqnn
\ar\ar\beta=\frac{1}{2}\I_{\{\sigma_2>0\}}+\frac{\theta_2-1}{2}\I_{\{\sigma_2=0, \int_0^x\xi^2\,n_2(\d \xi)\ge C_2x^{2-\theta_2}\I_{\{x\le1\}}\}}\in(0,1/2],\cr
\ar\ar\,\cr
\ar\ar
c_0=2^{-1}b_2\beta l_0^{\beta-1}l_1^{\alpha_2}\mrm{e}^{-(1-\beta)(l_1-l_0)/l_0},\quad c_1=(1-\beta)l_1/(\alpha_2l_0)
\eeqnn
and the constants $l_0,l_1$ will be determined later.  The following lemma is fundamental.

\begin{lem1}\label{le:3.5} The function $f$ defined by \eqref{f0} satisfies $f\in C^2(\mbb{R}_+)$.
  The derivatives $f'$ and $f''$ are given, respectively, by
\beqnn
f'(z) =
\begin{cases}
	\beta z^{\beta-1},&\quad z\in[0,l_0],\\
	\beta l_0^{\beta-1}\mrm{e}^{-(1-\beta)(z-l_0)/l_0}, &\quad z\in(l_0,l_1],\\
\frac{2c_0}{b_2}[c_1(z-l_1)+l_1]^{-\alpha_2}, &\quad z\in(l_1,\infty),
\end{cases}
\eeqnn
and
\beqnn
f''(z) =
\begin{cases}
	-\beta(1-\beta) z^{\beta-2},&\quad  z\in[0,l_0],\\
	-\beta(1-\beta)l_0^{\beta-2}\mrm{e}^{-(1-\beta)(z-l_0)/l_0}, &\quad z\in(l_0,l_1],\\
\frac{-2c_0c_1\alpha_2}{b_2}[c_1(z-l_1)+l_1]^{-\alpha_2-1}, &\quad z\in(l_1,\infty).
\end{cases}
\eeqnn
In particular,  $f'(z)\ge0$ and $f''(z)\le0$ for all $z\ge0$, and on each of the three intervals $(0,l_0)$, $(l_0, l_1)$ and $(l_1,\infty)$, we also have $f'''\ge0$ and $f^{(4)}\le0$. Moreover,
\beqnn
\sup_{z>0}f(z)=f(l_1)
+\frac{\beta\alpha_2}{(1-\beta)(\alpha_2-1)}l_0^{\beta}\mrm{e}^{-(1-\beta)(l_1-l_0)/l_0}<\infty.
\eeqnn
\end{lem1}

Given the function $f$ defined by \eqref{f0}, we can define a sequence of functions $\{F_n\}_{n\ge1} \subset {\mathcal{D}(L_{\bar{Z}})}$ such that
\begin{equation}\label{e:Funtion}
F_n(z)
\begin{cases}
	\le f(z) ,&\quad 0 < z\le 1/n,\\
    =   f(z), &\quad 1/n < z <  \infty.
\end{cases}
\end{equation}

\begin{pro1}\label{prop:3.6}
Suppose that $\alpha_2>1$ and either $\sigma_2>0$ or there are constants $\theta_2\in(1,2)$ and $C_2>0$ such that
$$
\int_0^x\xi^2\,n_2(\d \xi)\ge C_2x^{2-\theta_2},\quad x\in (0,1].
$$
Then, given any $M>0$, there exist constants $n_0>1$ and $C_3>0$ such that for any $n>n_0$,
\beqlb\label{ineq:3.12}
L_{\bar{Z}}F_n(z)\le-C_3\I_{\{z>1/n\}}.
\eeqlb
Consequently, for any $M>0$ and $t_0>2\sup_{z>0}f(z)/C_3$,
\beqlb\label{ineq:3.13}
\sup_{y>\tilde{y}\ge0}\mbb{P}_{(y, \tilde{y})}(\bar{\zeta}_0 \ge t_0)\le \frac{1}{2}.
\eeqlb
\end{pro1}

\proof
Let us recall that $$
\beta=\frac{1}{2}\I_{\{\sigma_2>0\}}+\frac{\theta_2-1}{2}\I_{\{\sigma_2=0, \int_0^x\xi^2\,n_2(\d \xi)\ge C_2x^{2-\theta_2}\I_{\{x\le1\}}\}}.
$$
Under the assumptions, we can choose a constant $\kappa_0>0$ such that
\begin{equation}\label{ineq:3.14}
\sigma_2+\frac{1}{3}\int_0^x\xi^2\,n_2(\d \xi)\ge \kappa_0x^{1-2\beta},\quad x\in (0,1].
\end{equation}
From \eqref{generator of bar Z} and the construction of $\{F_n\}_{n\ge1}$, for any $n\ge1$ and $z>1/n$, we have
\beqlb\label{ineq:3.15}\begin{split}
\bar{L}_{Z}F_n(z) \le&[(2kM+a^+_2)z-b_2z^{\alpha_2}]F_n'(z)+z\sigma_2 F_n''(z)
+z\int_0^\infty D_\xi F_n(z)\,n_2(\d \xi)\\
 = & [(2kM+a^+_2)z-b_2z^{\alpha_2}]f'(z)+z\sigma_2 f''(z)
+z\int_0^\infty D_\xi f(z)\,n_2(\d \xi).\end{split}
\eeqlb
Let $l_0\in(0,1-\beta)$,  $l_1>2$ and $n_0>1/l_0$ to be specified later.
Since $l_0<1$ and $l_1\ge2$, for any $n>n_0$, we divide $(1/n,\infty)$ into three intervals: $(1/n,l_0]$, $(l_0,l_1/2)$, $[l_1/2, \infty)$, and consider the corresponding estimate on each interval.

(i) We first deal with the case that $z\in(1/n,l_0]$.
Observe that $D_\xi f\le 0$, $f''\le0$ and $f'''\ge0$ 
 on $(0,l_1)$. Then, for any $z\in(1/n,l_0]$, combining \eqref{ineq:3.14}, \eqref{ineq:3.15} and the mean-value theorem gives that
\beqnn
\begin{split}
L_{\bar{Z}}F_n(z)
 \le &
 (2kM+a^+_2)zf'(z)
+z\sigma_2 f''(z)+z\int_0^z D_\xi f(z)\,n_2(\d \xi)\\
\le&(2kM+a^+_2)zf'(z)
+z\sigma_2 f''(2z)+\frac{1}{2}zf''(2z)\int_0^{z}\xi^2\,n_2(\d \xi)\\
 \le&(2kM+a^+_2)zf'(z)
+\kappa_0z^{2-2\beta}f''(2z).\end{split}
\eeqnn
It follows from  Lemma \ref{le:3.5} that
\beqnn
\begin{split} &L_{\bar{Z}}F_n(z)\\
&\le
\beta(2kM+a^+_2)z^\beta-\kappa_0\beta(1-\beta)z^{2-2\beta}
\left[(2z)^{\beta-2}\I_{\{z\le l_0/2\}}
+l_0^{\beta-2}\e^{(\beta-1)(2z-l_0)/l_0}\I_{\{l_0/2<z\le l_0\}}
\right]\\
&\le
\beta z^{-\beta}\left[
(2kM+a_2^+)l_0^{2\beta}-\kappa_0(1-\beta)\left(
2^{\beta-2}\I_{\{z\le l_0/2\}}+
\e^{\beta-1}2^{\beta-2}\I_{\{l_0/2<z\le l_0\}}
\right)
\right]\\
&\le
\beta z^{-\beta}\left[
(2kM+a_2^+)l_0^{2\beta}-\kappa_0(1-\beta)2^{\beta-2}\e^{\beta-1}
\right].\end{split}
\eeqnn
By choosing a sufficiently small  $l_0\in (0, 1-\beta)$ such that
$$
(2kM+a_2^+)l_0^{2\beta}\le \kappa_0(1-\beta)2^{\beta-3}\e^{\beta-1},
$$
we obtain
$$
L_{\bar{Z}}F_n(z)\le-\kappa_0\beta(1-\beta)2^{\beta-3}\e^{\beta-1}z^{-\beta}\le -\kappa_0\beta(1-\beta)2^{\beta-3}\e^{\beta-1}l_0^{-\beta}.
$$

(ii) Next we turn to the case that $z\in(l_0,l_1/2)$. Note that $l_0\le 1-\beta$ and $l_1\ge2$, we then have
\beqnn
z+l_0/(1-\beta)\le l_1/2+1< l_1.
\eeqnn
Since  $f^{(4)}\le0$ on $(l_0,l_1)$, the Taylor expansion implies that for any $\xi\in (0, l_0/(1-\beta))$,
\begin{align}\label{ineq:3.16}
D_\xi f(z)
\le\frac{1}{2}\xi^2 f''(z)+\frac{1}{6}\xi^3 f'''(z).
\end{align}

Using \eqref{ineq:3.16}, together with the fact that $D_{\xi}f(z)\le0$ for $z\in(l_0,l_1)$, we have
\beqnn
\int_0^\infty D_\xi f(z)\,n_2(\d \xi)
\ar\le\ar
\int_0^{l_0/(1-\beta)}\left(\frac{1}{2}\xi^2 f''(z)+\frac{1}{6}\xi^3 f'''(z)\right)\,n_2(\d \xi)\cr
\ar\le\ar
\int_0^{l_0/(1-\beta)}
\left(\frac{1}{2}\xi^2f''(z)+\frac{l_0}{6(1-\beta)}\xi^2f'''(z)\right)\,n_2(\d \xi)\cr
\ar\le\ar
\frac{1}{3}f''(z)\int_0^{l_0/(1-\beta)}\xi^2\,n_2(\d \xi),
\eeqnn
where in the  third inequality  we have used $f'''(z)=-(1-\beta)l_0^{-1}f''(z)\ge0$ for $z\in(l_0,l_1)$. Hence, combining \eqref{ineq:3.14}, \eqref{ineq:3.15}, Lemma  \ref{le:3.5} and the fact that $l_0/(1-\beta)\le1$,
\beqnn
L_{\bar{Z}}F_n(z)
\ar\ar\le
(2kM+a^+_2)zf'(z)+z\sigma_2f''(z)+\frac{1}{3}zf''(z)
\int_0^{l_0/(1-\beta)}
\xi^2\,n_2(\d \xi)\cr
\ar\ar\le
(2kM+a^+_2)zf'(z)+\kappa_0[l_0/(1-\beta)]^{1-2\beta}zf''(z)\cr
\ar\ar=
\beta l_0^{\beta-1}
\left(
2kM+a^+_2-\kappa_0[l_0/(1-\beta)]^{-2\beta}
\right)\mrm{e}^{-(1-\beta)(z-l_0)/l_0}z.
\eeqnn
By choosing a sufficiently small  $l_0\in (0, 1-\beta)$ such that
$$
2kM+a^+_2-\kappa_0[l_0/(1-\beta)]^{-2\beta}\le-1,
$$
we have
$$
L_{\bar{Z}}F_n(z)
\le-\beta l_0^{\beta-1}
\mrm{e}^{-(1-\beta)(z-l_0)/l_0}z
\le-\beta l_0^{\beta}
\mrm{e}^{-(1-\beta)(l_1/2-l_0)/l_0}.
$$

(iii) We finally consider the case that $z\in[l_1/2, \infty)$.  In view of \eqref{ineq:3.15}, it is obvious from the mean-value theorem and $f''\le0$  that
\beqlb\label{ineq:3.18}
L_{\bar{Z}}F_n(z)\le [(2kM+a^+_2)z-b_2z^{\alpha_2}]f'(z).
\eeqlb
By choosing a large enough $l_1\ge2$ such that for all $z\ge l_1/2$,
$$
(2kM+a^+_2)z-b_2z^{\alpha_2}\le-b_2z^{\alpha_2}/2,
$$
we have from \eqref{ineq:3.18} and Lemma \ref{le:3.5} that
\beqnn
\begin{split}
L_{\bar{Z}}F_n(z)&
\le-\frac{b_2}{2}z^{\alpha_2}f'(z)\\
&=-\frac{b_2
\beta l_0^{\beta-1}}{2}\mrm{e}^{-(1-\beta)(z-l_0)/l_0}z^{\alpha_2}\I_{\{z\in[l_1/2,l_1)\}}-
\frac{c_0z^{\alpha_2}}{[c_1(z-l_1)+l_1]^{\alpha_2}}\I_{\{z\in[l_1,\infty)\}}.\end{split}
\eeqnn
Thanks to the fact that
$$
\lim_{z\rightarrow\infty}
\frac{z^{\alpha_2}}{[c_1(z-l_1)+l_1]^{\alpha_2}}=c_1^{-\alpha_2},
$$
there exists a constant $c_2>0$ such that for any $z\ge l_1/2$,
$L_{\bar{Z}}F_n(z) \le -c_2.$

Combining all the estimates in the three cases above, we obtain \eqref{ineq:3.12} for some $C_3>0$.

We now prove \eqref{ineq:3.13}. By using It\^{o}'s formula and  \eqref{ineq:3.12}, for any $n>n_0$ and $y-\tilde{y} > 1/n$, we have
$$
\mbb{E}_{(y,\tilde{y})}\left[
F_n(\bar{Z}_{t\wedge\bar{\zeta}_{1/n}})\right] = F_n(y-\tilde{y})+\mbb{E}_{(y,\tilde{y})}\int_0^{t\wedge \bar{\zeta}_{1/n}}L_{\bar{Z}}F_n(\bar{Z}_s)\,\d s \le  F_n(y-\tilde{y})-C_3\mbb{E}_{(y,\tilde{y})}[t\wedge \bar{\zeta}_{1/n}].
$$
Letting $t \rightarrow \infty$, 
we deduce that
$$
\sup_{y-\tilde{y} >1/n}\mbb{E}_{(y,\tilde{y})}[\bar{\zeta}_{1/n}] \le \sup_{z>1/n}F_n(z)/C_3\le \sup_{z>0}f(z)/C_3.
$$
Therefore, by letting $n\rightarrow\infty$ and applying the Markov inequality, for any $t_0>2\sup_{z>0}f(z)/C_3$,  we have
$$
\sup_{y>\tilde{y}\ge 0}\mathbb{P}_{(y,\tilde{y})}(\bar{\zeta}_0\ge t_0)\le \frac{1}{2}.
$$
The proof is complete. \qed

The two local comparison principles in Section 3.2.1 are only valid, while $(X_t)_{t\ge0}$ remains below $2M$.  For this reason, we next estimate the exit time $\tau^+_{2M}$. Recall that $\mbb{P}_x(\cdot)$ is the law of the process $(X_t)_{t\ge0}$ starting from $x$.

\bglemma\label{le:3.7}
Given any $M > 1$, there is  a constant $t'>0$ such that
\beqnn
\mbb{P}_M(\tau^+_{3M/2} > t')>0.
\eeqnn
\edlemma

\proof Let $(X_t)_{t\ge0}$ be the solution to 
the first component of the SDE 
\eqref{main sde} starting from $X_0=M$. For $t>0$, define
$$
L_t=\sqrt{2\sigma_1}
\int_0^t\int_0^{X_s}\,W_1(\d s,\d u)+
\int_0^t\int_0^1\int_0^{X_{s-}}\xi\,\tilde{N}_1(\d s, \d \xi, \d u).
$$
From \eqref{main sde} we have
$$
X_{t\wedge\tau^+_{3M/2}}-M=L_{t\wedge\tau^+_{3M/2}}+X^{(1)}_{t\wedge\tau^+_{3M/2}}+X^{(2)}_{t\wedge\tau^+_{3M/2}},\quad t\ge0,
$$
where
$$
X^{(1)}_t=
\int_0^t
\left(-b_1X_s^{\alpha_1}+\left(a_1-\int_1^\infty\xi\,n_1(\d\xi)\right)X_s+\gamma_1\right)\,\d s
$$
and
$$
X^{(2)}_t=\int_0^t\int_1^\infty\int_0^{X_{s-}}\xi\,N_1(\d s, \d \xi, \d u).
$$
Then
\beqlb\label{ineq:3.19}
\begin{split} \mbb{P}_{M}(\tau^+_{3M/2} \le t)
&\le \mbb{P}_{M}(X_{\tau^+_{3M/2}\wedge t}-M\ge M/2)\\
&\le
\mbb{P}_M(L_{\tau^+_{3M/2}\wedge t}>M/4)
+\mbb{P}_M(X^{(1)}_{\tau^+_{3M/2}\wedge t}>M/4)+\mbb{P}_M(X^{(2)}_{\tau^+_{3M/2}\wedge t}>0).\end{split}
\eeqlb
Clearly, $(L^*_t)_{t\ge0}:=(L_{t\wedge\tau^+_{3M/2}})_{t\ge0}$ is a martingale with quadratic variation process
\beqnn
\langle L^*\rangle_t=\int_0^{t\wedge\tau^+_{3M/2}}
2\sigma_1X_s\,\d s
+\int_0^{t\wedge\tau^+_{3M/2}}X_s\,\d s\int_0^1\xi^2\,n_1(\d\xi),\quad t\ge0.
\eeqnn
Using Jensen's inequality and the Markov inequality, we have
\beqlb\label{ineq:3.20}
\begin{split}
\mbb{P}_M(|L^*_t|>M/4)  \le& \frac{4}{M}\mbb{E}_M|L^*_t| \le  \frac{4}{M}(\mbb{E}_M|L^*_t|^2)^{1/2} = \frac{4}{M}(\mbb{E}_M\langle L^* \rangle_t)^{1/2}\\
 =&\frac{4}{M}\Bigg\{
\mathbb{E}_M\left[\int_0^{t\wedge\tau^+_{3M/2}}
2\sigma_1X_s\,\d s\right]\\
&\quad\quad
+\mathbb{E}_M\left[
\int_0^{t\wedge\tau^+_{3M/2}}X_s\,\d s\int_0^1\xi^2\,n_1(\d\xi)
\right]
\Bigg\}^{1/2}\\
 \le&
4\left(3\sigma_1 +\frac{3}{2}\int_0^1\xi^2\,n_1(\d \xi)\right)^{1/2}t^{1/2}M^{-1/2}.\end{split}
\eeqlb

A direct calculation shows that
$$
X^{(1)}_{\tau^+_{3M/2}\wedge t}\le (3M|a_1|/2+\gamma_1)t
$$
and hence
\beqlb\label{ineq:3.21}
\mbb{P}_M(X^{(1)}_{\tau^+_{3M/2}\wedge t}>M/4)
\le 4\mbb{E}_{M}(X^{(1)}_{\tau^+_{3M/2}\wedge t})/M
\le 6|a_1|t+4\gamma_1t/M.
\eeqlb

Noting that $(X^{(2)}_t)_{t\ge0}$ is a pure jump process and that the waiting time before its first jump
can be controlled by an exponential random variable, we have
\beqlb\label{ineq:3.22}
\mbb{P}_M(X^{(2)}_{\tau^+_{3M/2}\wedge t}\neq0)
\le 1-\exp(-3Mn_1(1,\infty)t/2) \le 3n_1(1,\infty)tM/2.
\eeqlb

Therefore, by substituting \eqref{ineq:3.20}, \eqref{ineq:3.21} and \eqref{ineq:3.22} into \eqref{ineq:3.19}, there exists a constant $\lambda>0$  such that for any $t\in(0,1)$,
\beqnn
\mbb{P}_M(\tau^+_{3M/2}\le t)\le\lambda M t^{1/2}.
\eeqnn
Finally, for any $
0<t'< 1\wedge (\lambda M)^{-2},$
we obtain
$$
\mbb{P}_M(\tau^+_{3M/2}\le t')<1,
$$
which yields the  desired result.
\qed

Set
$$
\phi(x)=-b_1x^{\alpha_1}+a_1x+\gamma_1,\quad x\ge0.
$$
If $\alpha_1 > 1$, one can check the equation $\phi(x)=0$ has the unique positive solution $x_0$ so that $\phi(x)<0$ for all $x > x_0$. Recall that for any $b\ge0$, $\tau^-_b=\inf\{t>0: X_t\le b\}$.

\bglemma\label{le:3.8}
Suppose that $\alpha_1>1$. Given any $M>x_0$,
\beqnn
\mbb{P}_{3M/2}(\tau^-_{M}<\tau^+_{2M})>0.
\eeqnn
\edlemma
\proof We define a function
$$
w(x)=1-\frac{M^{-1}x}{1+M^{-1}x},\quad x\ge0.
$$
From this definition, it is obvious that $w(M) = \frac{1}{2}$, $w(3M/2)=\frac{2}{5}$, $w(2M) = \frac{1}{3}$, $w' < 0$ and $w'' >0$ on $\R_+$. Since $M>x_0$, we have $\phi(x)<0$ for all $x\ge M$. 
Combining the mean-value theorem with  \eqref{generator}
yields
\beqlb\label{ineq:3.23}
L_Xw(x)=\phi(x)w'(x) > 0, \quad x\in[M,2M].
\eeqlb
Now, applying Proposition \ref{prop:2.1} and using \eqref{ineq:3.23}, we obtain
\beqnn
\mbb{E}_{3M/2}[w(X_{\tau^-_{M}\wedge\tau^+_{2M}})] = w(3M/2)+\mbb{E}_{3M/2}\left[\int_0^{\tau^-_{M}\wedge\tau^+_{2M}}L_Xw(X_s)\,\d s\right]
\ge w(3M/2)=2/5.
\eeqnn
Therefore,
$$
 \mbb{P}_{3M/2}(\tau^-_M<\tau^+_{2M})+\frac{1}{3}\mbb{P}_{3M/2}(\tau^-_M\ge\tau^+_{2M})\ge \mbb{E}_{3M/2}[w(X_{\tau^-_{M}\wedge\tau^+_{2M}})] \ge 2/5,
$$ which implies that
$$ \mbb{P}_{3M/2}(\tau^-_M<\tau^+_{2M})\ge (2/5-1/3)/(1-1/3)=1/10.$$
The proof is complete.\qed

The estimates of Lemma \ref{le:3.7} and Lemma \ref{le:3.8} can now be combined through the strong Markov property, 
which in turn yields 
the tail estimate for $\tau^+_{2M}$.

\bgproposition\label{prop:3.9}
Suppose that $\alpha_1>1$. Then for any $M > (x_0 \vee 1)$ and $t>0$,
\beqnn
\mbb{P}_M(\tau^+_{2M}>t)>0.
\eeqnn
\edproposition

\proof Let $t'$ be the constant given in  Lemma \ref{le:3.7}. Using the strong Markov property at $\tau^{-}_M$, together with Lemma \ref{le:3.7} and Lemma \ref{le:3.8}, we see that
\beqnn
\mbb{P}_{3M/2}(\tau^+_{2M}>t')\ar\ge\ar\mbb{P}_{3M/2}(\tau^-_{M}<\tau^+_{2M}, \tau^+_{2M}\circ\theta_{\tau^-_M}>t')\cr
\ar=\ar
 \mbb{P}_{3M/2}(\tau^-_{M}<\tau^+_{2M})
 \mbb{P}_M(\tau^+_{2M}>t')\cr
 \ar\ge\ar
 \mbb{P}_{3M/2}(\tau^-_{M}<\tau^+_{2M})\mbb{P}_M(\tau^+_{3M/2}>t')>0,
\eeqnn
where in the equality above we used the strong Markov property at $\tau_M^-$ and the fact that the process $(X_t)_{t\ge0}$ has nonnegative jumps.
By Dawson and Li \cite[Theorem 2.3]{DaL12}, we have
$$
\mbb{P}_{M}(\tau^+_{2M}>t')>\mbb{P}_{3M/2}(\tau^+_{2M}>t')>0.
$$

We now extend this to arbitrary multiples of $t'$. Applying the Markov property at $t'$, we can establish that
\beqnn
\mbb{P}_M(\tau^+_{2M}>2t')\ar=\ar
\mbb{E}_M\left[\I_{\{\tau^+_{2M}>t'\}}\mbb{P}_{X_{t'}}(\tau^+_{2M}>t')\right]\cr
\ar\ge\ar
\mbb{E}_M\left[\I_{\{\tau^+_{3M/2}>t'\}}\mbb{P}_{X_{t'}}(\tau^+_{2M}>t')\right]\cr
\ar\ge\ar\mbb{E}_M\left[\I_{\{\tau^+_{3M/2}>t'\}}\mbb{P}_{3M/2}(\tau^+_{2M}>t')\right]>0,
\eeqnn where in the second inequality we used Dawson and Li \cite[Theorem 2.3]{DaL12} as well.
Hence, by induction, $$
\mbb{P}_M(\tau^+_{2M}>nt')>0,\quad  n=1,2,....
$$
Note that for any fixed $t>0$, there exists an integer $n_1\ge1$ such that
$$
(n_1-1)t'<t\le n_1t'.
$$
Then, for such $t$ we have
$$
\mbb{P}_M(\tau^+_{2M}>t)>\mbb{P}_M(\tau^+_{2M}>n_1t')>0.
$$
This completes the proof.\qed

\subsubsection{{\bf Proof of Proposition $\ref{prop:1.7}$}}\label{Section3.2.3}\noindent {\it Proof of Proposition $\ref{prop:1.7}$.} By Lemma \ref{le:3.2}, we only need to consider the case that $y>\tilde{y} \ge 0$. Set
 $Y_0=y$ and $\tilde Y_0=\tilde y$. Observe that $\bar{\zeta}_0$ and $\tau^+_{2M}$ are  independent, since the random elements $W_1, W_2$, $N_1$, and $N_2$ involved in are mutually independent. Then, for any $M> (x_0 \vee 1)$ and $t>0$, we have
\beqnn
\inf_{x\in[0,M],y>\tilde{y}\ge0}\mbb{P}_{(x,x,y,\tilde{y})}
(\bar{\zeta}_0 < t \le \tau^+_{2M}) \ar=\ar
\inf_{x\in[0,M],y>\tilde{y}\ge0}\left[
\mbb{P}_x
(\tau^+_{2M} \ge t)\mbb{P}_{(y,\tilde{y})}(\bar{\zeta}_0 < t)\right]\cr
\ar\ge\ar
\inf_{x\in[0,M]}\mbb{P}_x
(\tau^+_{2M} \ge t) \inf_{y>\tilde{y}\ge0}\mbb{P}_{(y,\tilde{y})}(\bar{\zeta}_0 < t).
\eeqnn
Recall that, by Dawson and Li \cite[Theorem 2.3]{DaL12}, the function $x\mapsto \mbb{P}_x
(\tau^+_{2M}\ge t)$ 
on $[0,M]$ attains  its minimum at $x=M$. This, together with Proposition \ref{prop:3.6} and Proposition \ref{prop:3.9}, implies that
\beqlb\label{ineq:3.24}
\inf_{x\in[0,M],y>\tilde{y}\ge0}\mbb{P}_{(x,x,y,\tilde{y})}
(\bar{\zeta}_0 < t_0 \le \tau^+_{2M}) \ge
\mbb{P}_M
(\tau^+_{2M}\ge t_0) \inf_{y>\tilde{y}\ge0}\mbb{P}_{(y,\tilde{y})}(\bar{\zeta}_0< t_0)>0,
\eeqlb where $t_0$ is the constant given in \eqref{ineq:3.13}.
On the other hand, according to Lemma \ref{le:3.3} and Lemma \ref{le:3.4}, we have
\beqlb\label{ineq:3.25}
\inf_{x\in[0,M], y>\tilde{y}\ge0}\mbb{P}_{(x,x,y,\tilde{y})}(\tilde{T}\wedge \tau_{2M}^+ \le \bar{\zeta}_0)=1.
\eeqlb
Combining \eqref{ineq:3.24} 
with  \eqref{ineq:3.25}, we therefore have
$$
\inf_{x\in[0,M],y>\tilde{y}\ge0}\mbb{P}_{(x,x,y,\tilde{y})}(\tilde{T}<\tau^+_{2M}\wedge t_0)
\ge
\inf_{x\in[0,M],y>\tilde{y}\ge0}\mbb{P}_{(x,x,y,\tilde{y})}
(\tilde{T}\wedge\tau^+_{2M} < t_0 \le \tau^+_{2M})>0.
$$
The proof is complete.\qed

\subsection{Proof of Proposition \ref{prop:1.8}}

\noindent {\it Proof of Proposition $\ref{prop:1.8}$.}  Let $(X_t, \tilde{X}_t)_{t \ge 0}$ be the coupling process of the first component process of the CBIPC-process given by \eqref{main sde}. By Dawson and Li \cite[Theorem 2.3]{DaL12}, $X_t\ge\tilde{X}_t$ for any $t\ge0$ almost surely, if $X_0\ge\tilde{X}_0 \ge 0$. Then, it suffices to consider the case $X_0 > \tilde{X}_0 \ge 0$. 
In particular, in this case, the coupling process  $(X_t, \tilde{X}_t)_{t \ge 0}$ has state space $D:=\{(x,\tilde{x}): x \ge \tilde{x}\ge0\}$ with generator $(\tilde{L}_X, {\mathcal{D}}(\tilde{L}_X))$. 
Here, for a function $F$ defined on $D$, being twice continuously differentiable on $D_0:=\{(x,\tilde{x}): x > \tilde{x}\ge0\}$, $\tilde{L}_XF(x,\tilde x)$ takes the form
\beqlb\label{coupling gene of X}
\begin{split}
\tilde{L}_X F(x, \tilde{x})
 =&(-b_1x^{\alpha_1}+a_1x+\gamma_1)F'_x(x, \tilde{x})
+(-b_1\tilde{x}^{\alpha_1}+a_1\tilde{x}+\gamma_1)F'_{\tilde{x}}(x, \tilde{x})\\
&
+\sigma_1xF''_{xx}(x, \tilde{x})
+\sigma_1\tilde{x}F''_{\tilde{x}\tilde{x}}(x, \tilde{x})
+2\sigma_1\tilde{x}F''_{x\tilde{x}}(x, \tilde{x})\\
&
+\tilde{x}\int_0^\infty
D_{(\xi,\xi)} F(x,\tilde{x})\,n_1(\d \xi)
 +(x-\tilde{x})\int_0^\infty D_{(\xi,0)} F(x,\tilde{x})\,n_1(\d \xi),\end{split}
 \eeqlb
 and ${\mathcal{D}}(\tilde{L}_X)$ denotes the linear space consisting of the functions $F$ such that the integrals above converge and define functions on locally bounded  on compact subsets of $D_0$.

Under Condition \ref{Condition}, we define the function
\beqnn
h(r) =
\begin{cases}
	1+r^\beta,&\quad r\in[0,l_0],\\
	h(l_0)+\frac{\beta l_0^\beta}{1-\beta}\left[1-\mrm{e}^{-\frac{1-\beta}{l_0}(r-l_0)}\right], & \quad r\in(l_0,l_1],\\
h(l_1)+c_0\int_0^{r-l_1}\frac{2}{b_1(c_1s+l_1)^{\alpha_1}}\,\d s, &\quad r\in(l_1,\infty).
\end{cases}
\eeqnn
 Here the parameters involved in the definition of $h$ are identical to those in \eqref{f0}, except that we replace $(b_2,\alpha_2,a_2,\gamma_2,\sigma_2,n_2)$ with  $(b_1,\alpha_1,a_1,\gamma_1,\sigma_1,n_1)$. Thanks to $\alpha_1>1$, we have
 \beqlb\label{ineq:3.27}
x^{\alpha_1}-\tilde{x}^{\alpha_1} \ge (x-\tilde{x}) x^{\alpha_1-1} \ge (x-\tilde{x})^{\alpha_1},\quad (x,\tilde{x}) \in D_0.
 \eeqlb
Let
 $\{H_n\}_{n\ge1}
 \subset {\mathcal{D}}(\tilde{L}_X)$ be a sequence of functions satisfying
 \beqnn
H_n(x,\tilde{x})
\begin{cases}
	\le h(x-\tilde{x}),&\quad 0 < x-\tilde{x} \le 1/n,\\
	= h(x-\tilde{x}), &\quad 1/n < x-\tilde{x} <\infty.
\end{cases}
\eeqnn
This, together with \eqref{coupling gene of X} and \eqref{ineq:3.27}, implies that, for $x-\tilde{x}>1/n$,
\beqnn
\tilde{L}_XH_n(x,\tilde{x})
\ar \le \ar\left[-b_1(x-\tilde{x})^{\alpha_1}+a_1(x-\tilde{x})\right]h'(x-\tilde{x})+\sigma_1(x-\tilde{x})h''(x-\tilde{x})\cr
\ar\ar
+(x-\tilde{x})\int_0^\infty
\left[h(x-\tilde{x}+\xi)-h(x-\tilde{x})-\xi h'(x-\tilde{x})\right]\,n_2(\d \xi).
\eeqnn
Following the proof of Proposition \ref{prop:3.6}, there exist constants $C_4>0$ and $n_0>1$
such that for all $n > n_0$,
 $$
\tilde{L}_XH_n(x,\tilde{x})
 \le-C_4\I_{\{x-\tilde{x}>1/n\}}.
 $$
 Moreover,  for any $t_0 > 2\sup_{r>0} h(r)/C_4$,
$$
\sup_{x>\tilde{x}\ge0}\mbb{P}_{(x,\tilde{x})}(\tilde{T}_X\ge t_0) < \frac{1}{2}.
$$
This completes the proof. \qed

\subsection{Proof of Remark \ref{remark on k le 0}}

Let $(X_t,Y_t)_{t\ge0}$ be the CBIPC-process satisfying \eqref{main sde} with $k\le0$, and let $(X_t, \tilde{X}_t, Y_t, \tilde{Y}_t)_{t\ge0}$ be the Markovian coupling process  with $X_0 \ge \tilde{X}_0 \ge0$.
Recall that
$$
\tilde{T}_X=\{t\ge0: X_t=\tilde{X}_t\},\quad  \tilde{T}=\inf\{t\ge0: (X_t,Y_t)=(\tilde{X}_t, \tilde{Y}_t)\}.
$$
For any $t \ge 0$,  it is easy to see that
\beqnn
\sup_{x,\tilde{x},y,\tilde{y}}\mbb{P}_{(x,\tilde{x},y,\tilde{y})}(\tilde{T} > 2t)  \leq \sup_{x,\tilde{x}}\mbb{P}_{(x,\tilde{x})}(\tilde{T}_X > t) + \sup_{x,\tilde{x},y,\tilde{y}}\mbb{P}_{(x,\tilde{x},y,\tilde{y})}( \tilde{T}_X \leq t, \,\, \tilde{T} > 2t).
\eeqnn
By the strong Markov property at $\tilde{T}_X$, we have
\beqnn
\sup_{x,\tilde{x},y,\tilde{y}}\mbb{P}_{(x,\tilde{x},y,\tilde{y})}( \tilde{T}_X \leq t, \,\, \tilde{T} > 2t) \ar \le \ar
\sup_{x,\tilde{x},y,\tilde{y}}\mbb{P}_{(x,\tilde{x},y,\tilde{y})}( \tilde{T}-\tilde{T}_X \ge t)\cr
\ar=\ar
\sup_{x,\tilde{x},y,\tilde{y}}\mbb{E}_{(x,\tilde{x},y,\tilde{y})}\left(\mbb{P}_{(x,\tilde{x},y,\tilde{y})}\left(\tilde{T}-\tilde{T}_X\ge t \mid\mathcal{F}_{\tilde{T}_X}\right)\right)\cr
\ar=\ar
\sup_{x,\tilde{x},y,\tilde{y}}\mbb{E}_{(x,\tilde{x},y,\tilde{y})}\left(\mbb{P}_{(X_{\tilde{T}_X},\tilde{X}_{\tilde{T}_X},Y_{\tilde{T}_X},\tilde{Y}_{\tilde{T}_X})}(\tilde{T}\ge t)\right)\cr
\ar\le\ar\sup_{x,y,\tilde{y}}\mbb{P}_{(x,x,y,\tilde{y})}(\tilde{T}\ge t),
\eeqnn
where in the last inequality we 
used the fact that $X_{\tilde{T}_X}=\tilde{X}_{\tilde{T}_X}$.
It follows that for any $t\ge0$,
\beqlb\label{ineq:3.28}
\sup_{x,\tilde{x},y,\tilde{y}}
\mbb{P}_{(x,\tilde{x},y,\tilde{y})}(\tilde{T} > 2t) \le \sup_{x,\tilde{x}}\mbb{P}_{(x,\tilde{x})}(\tilde{T}_X > t) + \sup_{x,y,\tilde{y}}\mbb{P}_{(x,x,y,\tilde{y})}(\tilde{T}\ge t).
\eeqlb

From the proof of Proposition \ref{prop:1.8}, for any $t_0>2\sup_{r>0}h(r)/C_4$, it holds that
\beqlb\label{ineq:3.29}
\sup_{x,\tilde{x}}\mbb{P}_{(x,\tilde{x})}(\tilde{T}_X \ge t_0) < \frac{1}{2}.
\eeqlb
On the other hand, let $\{F_n\}_{n\ge1}$ be the sequence of the functions defined by \eqref{e:Funtion}. Then, for any $x\ge0, n\ge1$ and $y-\tilde{y}>1/n$,
\beqlb\label{coupling gene of Y}\begin{split}
\tilde{L}_{x,x,Y}F_n(y,\tilde{y}) =&(kx+a_2)(y-\tilde{y})f'(y-\tilde{y})-b_2(y^{\alpha_2}-\tilde{y}^{\alpha_2})f'(y-\tilde{y})\\
&+\sigma_2(y-\tilde{y})f''(y-\tilde{y})
+(y-\tilde{y})\int_0^\infty
D_\xi f(y-\tilde{y})
\,n_2(\d \xi),\end{split}
\eeqlb
where $\tilde{L}_{x,x,Y}$ is defined in \eqref{coupling operator}. Since $\alpha_2>1$,  it follows that
$$
y^{\alpha_2}-\tilde{y}^{\alpha_2}\ge (y-\tilde{y})y^{\alpha_2-1}\ge(y-\tilde{y})^{\alpha_2}, \quad y\ge \tilde y\ge0.
$$
Combining this with \eqref{coupling gene of Y}, $k\le0$ and the fact that $f'>0$, we obtain that, for any $x\ge0, n\ge1$ and $y-\tilde{y}>1/n$,
\beqnn
\tilde{L}_{x,x,Y}G_n(y,\tilde{y})
\ar\le\ar
\left[a^+_{2}(y-\tilde{y})-b_2(y-\tilde{y})^{\alpha_2}\right]f'(y-\tilde{y})
+\sigma_2(y-\tilde{y})f''(y-\tilde{y})\cr
\ar\ar
+(y-\tilde{y})\int_0^\infty D_\xi f(y-\tilde{y})\,n_2(\d \xi).
\eeqnn
Then, following  the proof of Proposition \ref{prop:3.6}, there exist constants $C_5>0$ and $n_0>1$ such that for any $x\ge0$ and $n > n_0$,
$$
\tilde{L}_{x,x,Y}F_n(y,\tilde{y})\le-C_5\I_{\{y-\tilde{y}>1/n\}}.
$$
Moreover, for any $t_0> 2\sup_{r>0}f(r)/C_5$,
\beqlb\label{ineq:3.31}
\sup_{x,y, \tilde{y}}\mbb{P}_{(x,x,y,\tilde{y})}(\tilde{T}\ge t_0) < \frac{1}{2}.
\eeqlb

Finally, combining \eqref{ineq:3.28}, \eqref{ineq:3.29} and \eqref{ineq:3.31}, there exists some constant $\varepsilon\in(0,1)$
such that for any $t_0 >2(C_4\wedge C_5)^{-1}\sup_{r>0}(f(r)\vee h(r))$,
$$
\sup_{x,\tilde{x},y,\tilde{y}}\mbb{P}_
{(x,\tilde{x},y,\tilde{y})}(\tilde{T} > 2t_0) \le \varepsilon.
$$
The proof is then completed by following the proof of Theorem \ref{main result1}.
\qed

\section{Appendix}
 \setcounter{equation}{0}
 In this section, we  prove Proposition \ref{unique strong solution}, which gives the construction of the CBIPC-process $(X_t,Y_t)_{t\ge0}$. Due to the specific form of the interaction term, the existing results in the literature do not directly apply to our system. To overcome this, we employ a truncation argument. For any $N \ge 1$, let us consider
the following system of stochastic equations with any ${\mathcal{F}}_0$-measurable random vector $(X_0, Y_0) \in \mathbb{R}_+^2$:
\beqlb\label{trunca sde}
\left\{\begin{aligned}
X_t &= X_0 + \int_0^t \left(-b_1 X^{\alpha_1}_s + a_1 X_s + \gamma_1\right) \d s + \sqrt{2\sigma_1}\int_0^t\int_0  ^{X_s}\,W_1(\d s,\d u) \cr
&\quad + \int_0^t \int_0^\infty \int_0^{X_{s-}} \xi  \tilde{N}_1(\d s, \d \xi, \d u), \\[2mm]
Y^N_{t} &= Y_0 + \int_0^t \left[k (X_s\wedge N) Y^N_s - b_2 (Y^N_s)^{\alpha_2} + a_2 Y^N_s + \gamma_2\right] \mathrm{d}s \cr
&\quad + \sqrt{2\sigma_2}\int_0^t\int_0^{Y^N_s}\, W_2(\d s,\d u) + \int_0^t \int_0^\infty \int_0^{Y^N_{s-}} \xi  \tilde{N}_2(\d s, \d \xi, \d u),
\end{aligned}\right.
\eeqlb
where the parameters and the noises are the same as those in \eqref{main sde}.

Before proceeding to the proof, we first verify that the truncated system \eqref{trunca sde} preserves the non-negativity of solutions.

\bglemma
For any $N\ge1$, if $(X_t,Y^N_{t})_{t\ge0}$ satisfies \eqref{trunca sde}, then $\mathbb{P}(X_t\ge0,Y^N_t\ge0~~\text{for~all}~t\ge0)=1.$
\edlemma

\proof The assertion regarding $X:=(X_t)_{t\ge0}$ follows from the proof of Fu and Li \cite[Proposition 2.1]{FL10}. For any $N\ge1$, suppose that $Y^N:=(Y^N_t)_{t\ge0}$ leaves $[0,\infty)$ on some event $A$ with $\mbb{P}(A)>0$. Note that the jump term in the second component of the SDE \eqref{trunca sde} prevents $Y^N$ from jumping into $(-\infty,0)$. Then, on the event $A$, there exists a time interval $[t_0,t_1]$ such that $Y^N_{t_0}=0$ and $t\mapsto Y^N_t$ is a strictly negative continuous function on $[t_0,t_1]$.
 On the other hand, there are some $t'\in[t_0,t_1]$ and some $\varepsilon>0$ such that for all $s\in[t_0,t']$, $k (X_s\wedge N) Y^N_s-b_2(Y^N_s)^{\alpha_2}+a_2Y^N_s+\gamma_2 > \varepsilon$. Then, for all $t\in [t_0,t']$,
$$
Y^N_t=Y^N_t-Y^N_{t_0}=\int_{t_0}^t[k(X_s\wedge N)Y^N_s-b_2(Y^N_s)^{\alpha_2}+a_2Y^N_s+\gamma_2]\,\d s>\varepsilon(t-t_0)>0.
$$
Since $Y^N_t\le0$ for all $t\in[t_0,t_1]$, we get a contraction.\qed

\bgproposition\label{unique strong solution to local version}
For any $N\ge1$, there exists a unique nonnegative strong solution $(X_t,Y^N_t)_{t\ge0}$ to \eqref{trunca sde}.
\edproposition
\proof  By Dawson and Li \cite[Theorem 2.5]{DaL12}, there exists a unique nonnegative strong solution $(X_t)_{t\ge0}$ to the first component equation in \eqref{trunca sde}.

We now turn to the second component equation. First, following the proof of Lemma \ref{le:3.2}, one can show the pathwise uniqueness of the solution.
Then, it suffices to show that there exists a nonnegative weak solution to this equation. For any $n \ge 1$,  let $V_n=\{z\in\mathbb{R}_+: z \ge 1/n\}$. By \eqref{measure moment}, we have $n_1(V_n)+n_2(V_n)<\infty$. For $m,n\ge1$, by the results for continuous-type stochastic equations as in  Ikeda and Watanabe \cite[Theorem 2.2]{IW89}, one can show there is a weak solution to
\beqlb\label{local solution XY}
\left\{\begin{aligned}
 X_t&=X_0+\int_0^t\left(-b_1(X_s\wedge m)^{\alpha_1}+a_1(X_s\wedge m)+\gamma_1\right)\d s\cr
 &\quad + \sqrt{2{\sigma_1}}\int_0^t\int_0^{X_s\wedge m} W_1(\d s, \d u)\cr
&\quad-\left(\int_0^\infty(z-z\wedge m)\,n_1(\d z)+\int_{V_n}(z\wedge m)\,n_1(\d z)\right)\int_0^t(X_s\wedge m)\,\d s,\cr
Y^N_t&=Y_0+\int_0^t\left[k(X_s\wedge N\wedge m)(Y^N_s\wedge m)-b_2(Y^N_s\wedge m)^{\alpha_2}+a_2(Y^N_s\wedge m)+\gamma_2\right]\d s\cr
&\quad+\sqrt{2{\sigma_2}}\int_0^t\int_0^{Y^N_s\wedge m} W_2(\d s, \d u)\cr
&\quad-\left(\int_0^\infty(z-z\wedge m)\,n_2(\d z)+\int_{V_n}(z\wedge m)\,n_2(\d z)\right)\int_0^t(Y^N_s\wedge m)\,\d s.
\end{aligned}\right.
 \eeqlb
Again the pathwise uniqueness holds for the above system by the proof of Lemma \ref{le:3.2}. Then it has a unique strong solution to \eqref{local solution XY}.

Following the proof of Fu and Li \cite[Proposition 2.2]{FL10}, one can further show that there is a pathwise unique nonnegative strong solution $(X^{m,n}_t,Y^{N,m,n}_t)_{t\ge0}$ to
\beqlb\label{local solution XY1}
\left\{\begin{aligned}
 X_t &= X_0+\int_0^t\left(-b_1(X_s\wedge m)^{\alpha_1}+a_1(X_s\wedge m)+\gamma_1\right)\d s\cr
 &\quad+\sqrt{2\sigma_1}\int_0^t\int_0^{X_s\wedge m} W_1(\d s,\d u)-\int_0^\infty(z-z\wedge m)\,n_1(\d z)\int_0^t(X_s\wedge m)\,\d s\cr
&\quad +\int_0^t\int_{V_n}\int_0^{X_{s-}\wedge m}(z\wedge m)\,\tilde{N}_1(\d s,\d z,\d u),\cr
Y^N_t& = Y_0+\int_0^t\left[k(X_s\wedge N\wedge m)(Y^N_s\wedge m)-b_2(Y^N_s\wedge m)^{\alpha_2}+a_2(Y^N_s\wedge m)+\gamma_2\right]\d s\cr
&\quad+\sqrt{2\sigma_2}\int_0^t\int_0^{Y^N_s\wedge m} W_2(\d s, \d u)-\int_0^\infty(z-z\wedge m)\,n_2(\d z)\int_0^t(Y^N_s\wedge m)\,\d s\cr
 &\quad+\int_0^t\int_{V_n}\int_0^{Y^N_{s-}\wedge m}(z\wedge m)\,\tilde{N}_2(\d s,\d z,\d u).
\end{aligned}
\right.
 \eeqlb
As in the proof of Fu and Li \cite[Lemma 4.3]{FL10}, the sequence $\{(X^{m,n}_t,Y^{N,m,n}_t)_{t\ge0}: n=1,2,...\}$ is tight in $D([0,\infty),\mathbb{R}_+^2)$. Furthermore, following the proof of Fu and Li \cite[Theorem 4.4]{FL10}, it is easy to show that any weak limit point $(X^{m}_t,Y^{N,m}_t)_{t\ge0}$ of the sequence above is a nonnegative weak solution to
\beqlb\label{local solution XY2}
\left\{\begin{aligned}
 X_t &= X_0+\int_0^t\left(-b_1(X_s\wedge m)^{\alpha_1}+a_1(X_s\wedge m)+\gamma_1\right)\d s\cr
 &\quad+\sqrt{2\sigma_1}\int_0^t\int_0^{X_s\wedge m} W_1(\d s,\d u)-\int_0^\infty(z-z\wedge m)\,n_1(\d z)\int_0^t(X_s\wedge m)\,\d s\cr
& \quad+\int_0^t\int_0^\infty\int_0^{X_{s-}\wedge m}(z\wedge m)\,\tilde{N}_1(\d s,\d z,\d u),\cr
Y^N_t& = Y_0+\int_0^t\left[k(X_s\wedge N\wedge m)(Y^N_s\wedge m)-b_2(Y^N_s\wedge m)^{\alpha_2}+a_2(Y^N_s\wedge m)+\gamma_2\right]\d s\cr
&\quad+\sqrt{2\sigma_2}\int_0^t\int_0^{Y^N_s\wedge m} W_2(\d s, \d u)-\int_0^\infty(z-z\wedge m)\,n_2(\d z)\int_0^t(Y^N_s\wedge m)\,\d s\cr
 &\quad+\int_0^t\int_0^\infty\int_0^{Y^N_{s-}\wedge m}(z\wedge m)\,\tilde{N}_2(\d s,\d z,\d u).
\end{aligned}
\right.
 \eeqlb
According to the proof of Lemma \ref{le:3.2} again, the pathwise uniqueness holds for \eqref{local solution XY2}, so the system \eqref{local solution XY2} has a unique strong solution. Then, by a modification of the proof of Fu and Li \cite[Proposition 2.4]{FL10}, \eqref{trunca sde} has a unique strong solution.\qed

\noindent{\it Proof of Proposition $\ref{unique strong solution}$}.\,\,  For any $N\ge1$, by Proposition \ref{unique strong solution to local version}, there exists a unique nonnegative strong solution to \eqref{trunca sde}. For each $N\ge1$, recall that $\tau^+_N=\inf \{t>0: X_t\ge N\}$. By Fu and Li \cite[Proposition 2.3]{FL10}, we have $\tau^+_{N} \uparrow \infty$ almost surely as $N\rightarrow\infty$. Define $Y_t$ for $t < \tau^+_N$ by $Y_t = Y^N_t$. By the pathwise uniqueness of the solution to \eqref{trunca sde}, these definitions are consistent for different $N$; that is, for $M > N$, we have $Y_t^M = Y_t^N$ for all $t < \tau_N$ almost surely. Hence, we can define a system $(X_t, Y_t)_{t\ge0}$ by $Y_t=Y^N_t$ for $0\le t < \tau^+_N$, and this system satisfies \eqref{main sde} for all $0 \le t < \tau^+_N$. Since $\tau^+_N \uparrow \infty$ almost surely as $N\to\infty$, it follows that $(X_t, Y_t)_{t\ge 0}$ satisfies \eqref{main sde} for all $t\ge 0$ almost surely. The pathwise uniqueness follows from the uniqueness of the  system \eqref{trunca sde} and the fact that any solution to \eqref{main sde} must coincide with \eqref{trunca sde} on each stochastic interval $[0,\tau_N)$. This completes the proof. \qed

\ \

\noindent {\bf Acknowledgements.}\,\,
The research of Shukai Chen is supported by the National Natural Science Foundation of China (No.\ 12401167) and Fujian Provincial Natural Science Foundation of China (No.\ 2024J08050). The research of Pei-Sen Li is supported by the National Natural Science Foundation of China (No.\ 12271029). The research of Jian Wang is supported by the National Key R\&D Program of China (2022YFA1006003) and the National Natural Science Foundations of China (Nos.\ 12225104 and 12531007).

 \end{document}